\documentclass[11pt]{article}
\usepackage{amssymb,amsmath,amsfonts}
\usepackage{url}
\usepackage{cases}
\usepackage{hyperref}
\usepackage{xcolor}
\usepackage{graphicx}

\def\bT{{\mathbf{T}}}
\def\bS{{\mathbf{S}}}
\def\nuc{{\hbox{\tiny \rm nuc}}}

\newcommand{\qed}{\hfill $\Box$ \nr \medskip}
\def\mR{{{\mathfrak{R}}}}

\def\lin{{\hbox{\rm\tiny lin}}}

\def\s{{\hbox{\rm\scriptsize s}}}

\def\three?{3}
\def\four{\tfrac{1}{4}}
\def\ten?{10}
\def\bM{{\mathbf{M}}}
\def\bK{{\mathbf{K}}}

\def\Argmin{\mathop{\hbox{\rm Argmin}}}
\def\beq{\begin{equation}}
\def\eeq{\end{equation}}
\newtheorem{theorem}{Theorem}[section]

\newtheorem{proposition}{Proposition}[section]
\newtheorem{corollary}{Corollary}[section]

\def\norm2to2{{\|\cdot\|_{2,2}}}
\def\Prob{\hbox{\rm Prob}}

\def\bE{{\mathbf{E}}}

\def\Diag{\hbox{\rm  Diag}}
\def\Prob{\hbox{\rm  Prob}}

\def\Opt{\hbox{\rm Opt}}

\def\Conv{\hbox{\rm  Conv}}

\def\Tr{{\mathop{\hbox{\rm  Tr}}}}
\def\cA{{\cal A}}
\def\cB{{\cal B}}

\def\cH{{\cal H}}

\def\cL{{\cal L}}
\def\cM{{\cal M}}
\def\cN{{\cal N}}

\def\cQ{{\cal Q}}
\def\cR{{\cal R}}
\def\cS{{\cal S}}
\def\cT{{\cal T}}
\def\cU{{\cal U}}
\def\cV{{\cal V}}
\def\cW{{\cal W}}
\def\cX{{\cal X}}
\def\cY{{\cal Y}}
\def\cZ{{\cal Z}}

\def\rank{{\mathop{\hbox{\rm  Rank}}}}

\def\Argmin{\mathop{\hbox{\rm  Argmin}}}

\def\Det{{\mathop{\hbox{\rm  Det}}}}

\def\bK{{\mathbf{K}}}

\def\bS{{\mathbf{S}}}

\def\e{{\hbox{\rm e}}}

\def\qed{\ \hfill$\square$\par\smallskip}

\def\mypict3{\epsfxsize=220pt\epsfysize=80pt\epsfbox}

\def\bR{{\mathbf{R}}}

\newtheorem{lemma}{Lemma}[section]
\def\cH{{\cal H}}
\def\Col{{\hbox{\rm Col}}}

\def\Risk{{\hbox{\rm Risk}}}

\newcommand{\hide}[1]{{}}
\def\lin{{\mathrm{lin}}}
\def\poly{{\mathrm{poly}}}
\newcommand{\be}{\begin{eqnarray}}
\newcommand{\ee}[1]{\label{#1}\end{eqnarray}}
\newcommand{\nn}{\nonumber \\}
\newcommand{\ese}{\end{align*}}
\newcommand{\bse}{\begin{align*}}
\newcommand{\rf}[1]{~(\ref{#1})}
\newcommand{\wh}[1]{{\widehat{#1}}}
\def\mR{{\mathfrak{R}}}
\def\SG{{\cal SG}}
\def\half{\tfrac{1}{2}}
\def\wt#1{\widetilde{#1}}
\def\ov#1{\overline{#1}}

\def\mR{{{\mathfrak{R}}}}
\def\mr{{{\mathfrak{r}}}}
\def\ms{{{\mathfrak{s}}}}

\def\mP{{\mathfrak{p}}}

\def\Fro{\mathrm{Fro}}

\DeclareMathOperator*{\med}{median}

\title{On Robust Recovery of Signals from Indirect Observations}
\author{Yannis Bekri\thanks{\scriptsize LJK, Universit\'e Grenoble Alpes, Campus de Saint-Martin-d'H\`{e}res, 38401 France}
\and
Anatoli Juditsky$^*$
\and Arkadi Nemirovski
\thanks{\scriptsize Georgia Institute
 of Technology, Atlanta, Georgia
30332, USA}
}
\begin{document}
\maketitle
\begin{abstract}
Our focus is on robust recovery algorithms in statistical linear inverse problem. We consider two recovery routines---the much-studied linear estimate originating from Kuks and Olman \cite{KuksOlman} and polyhedral estimate introduced in \cite{juditsky2020polyhedral}. It was shown in \cite{PUP} that risk of these estimates can be tightly upper-bounded for a wide range of a priori information about the model through solving a convex optimization problem, leading to a computationally efficient implementation of nearly optimal estimates of these types. The subject of the present paper is design and analysis of linear and polyhedral estimates which are robust with respect to the uncertainty in the observation matrix. We evaluate performance of robust estimates under stochastic and deterministic matrix uncertainty and show how the estimation risk can be bounded by the optimal value of efficiently solvable convex optimization problem; ``presumably good'' estimates of both types are then obtained through optimization of the risk bounds with respect to estimate parameters.
\\~\\{\em 2020 Mathematics Subject Classification:} 62G05, 62G10, 90C90
\\{\em Keywords:} statistical linear inverse problems, robust estimation, observation matrix uncertainty
\end{abstract}

\section{Introduction}
In this paper we focus on the problem of recovering unknown signal $x$ given noisy observation $\omega\in \bR^m$,
\beq\label{eq01}
\omega=Ax+\xi,
\eeq
of the linear image $Ax$ of $x$; here $\xi\in \bR^m$  is observation noise. Our objective is to estimate the linear image $w=Bx\in\bR^\nu$ of
$x$ known to
belong to given convex and compact subset $\cX$ of $\bR^n$. The estimation problem above is a classical linear inverse problem.  When statistically analysed, popular approaches to solving \rf{eq01} (cf., e.g., \cite{natterer1986mathematics,johnstone1990speed,johnstone1991discretization,mair1996statistical,vogel2002computational,kaipio2006statistical,hoffmann2008nonlinear,proksch2018multiscale}) usually assume a special structure of the problem, when matrix $A$ and set $\cX$ ``fit each other,'' e.g., there exists a sparse approximation of the set $\cX$ in a given basis/pair of bases, in which matrix $A$ is ``almost diagonal'' (see, e.g. \cite{donoho1995nonlinear,cohen2004adaptive} for detail). Under these assumptions, traditional results focus on estimation algorithms which are both numerically straightforward and statistically (asymptotically) optimal with closed form analytical description of estimates and corresponding risks.
In this paper, $A$ and $B$ are ``general'' matrices of appropriate dimensions, and $\cX$ is a rather general convex and compact set.
Instead of deriving closed form expressions for estimates and risks (which under the circumstances seems to be impossible), we adopt an ``operational'' approach initiated in \cite{Don95} and further developed in \cite{JN2009,l2estimation,juditsky2020polyhedral,PUP}, within which
 both the estimate and its risk are yielded by efficient computation, rather than by an explicit analytical description.

In particular,  two classes of estimates were analyzed in \cite{l2estimation,juditsky2020polyhedral,PUP} in the operational framework.
\begin{itemize}
\item {\sl Linear estimates.} Since their introduction in \cite{kuks1,kuks2}, {\em linear estimates} are a standard part of the theoretical statistical toolkit. There is an extensive literature dealing with the design and performance analysis of linear estimates (see, e.g., \cite{Pinsker1980,donoho1990minimax,efromovich1996sharp,efromovich2008nonparametric,IbrHas1981,Tsybakov,wasserman2006all}). When applied in the estimation problem we consider here,
linear estimate $\widehat{w}_\lin^H(\omega)$ is of the form $\widehat{w}_H(\omega)=H^T\omega$
and is specified by a {\sl contrast matrix} $H\in\bR^{m\times \nu}$.
\item {\sl Polyhedral estimates.} The idea of a {\em polyhedral estimate} goes back to \cite{oldpaper} where it was
shown (see also \cite[Chapter 2]{saintflour}) that such estimate is near-optimal when recovering smooth multivariate regression function known to belong
to a given Sobolev ball from noisy observations taken along a regular grid. It has been recently reintroduced in \cite{grasmair2018variational}
and  \cite{proksch2018multiscale} and extended to the setting to follow in \cite{juditsky2020polyhedral}.
    In this setting, a polyhedral estimate $\omega\mapsto \widehat{w}^H_\poly(\omega)$ is specified by a {\sl contrast matrix} $H\in\bR^{m\times M}$ according to
$$
\omega\mapsto\wh x^H(\omega)\in\Argmin_{x\in\cX} \|H^T(\omega-Ax)\|_\infty\mapsto \widehat{w}^H_\poly(\omega) :=B\wh{x}(\omega).
$$
\end{itemize}
Our interest in these estimates stems from the results of \cite{JudNem2018,juditsky2020polyhedral,PUP} where it is shown that in the Gaussian case ($\xi\sim\cN(0,\sigma^2I_m)$), linear and polyhedral estimates with properly designed efficiently
computable contrast matrices are near-minimax optimal in terms of their risks over a rather general class of loss functions and signal sets---ellitopes
and spectratopes.
\footnote{Exact definitions of these sets are reproduced in the main body of the paper.
For the time being, it suffices to point out two instructive examples: the bounded
intersections of finitely many sets of the form $\{x:\|Px\|_p\leq1\}$, $p\geq2$, is an ellitope (and a spectratope as well),
and the unit ball  of the spectral norm in the space of $m\times n$ matrices is a spectratope.}

In this paper we consider an estimation problem which is a generalization of that mentioned above in which observation matrix
$A\in \bR^{m\times n}$ is {\em uncertain}. Specifically, we assume that
\beq\label{eq1}
\omega=A[\eta]x+\xi\eeq
where
$\xi\in\bR^m$ is zero mean {random noise} and
\beq\label{eq2}
A[\eta]=A+{\sum}_{\alpha=1}^q\eta_\alpha A_\alpha\in\bR^{m\times n}
\eeq
where $A,A_1,...,A_q$ are given matrices and $\eta\in \bR^q$ is uncertain perturbation (``uncertainty'' for short).
We consider separately two situations: the first one in which the perturbation $\eta$ is random (``random perturbation''),
and the second one with $\eta$ selected, perhaps in an adversarial fashion, from a given uncertainty set $\cU$ (''uncertain-but-bounded perturbation'').
Observation model \rf{eq1} with random uncertainty is related to the linear regression problem with random errors in regressors \cite{bennani1988utilisation,carroll1996use,fan1993nonparametric,gleser1981estimation,kukush2005consistency,stewart1990stochastic,van2013total} which is usually addressed through total least squares. It can also be seen as alternative modeling of the statistical inverse problem in which sensing matrix is recovered with stochastic error (see, e.g., \cite{cavalier2005adaptive,cavalier2007wavelet,efromovich2001inverse,hall2005nonparametric,hoffmann2008nonlinear,marteau2006regularization}).
 Estimation from observations \rf{eq1} under uncertain-but-bounded perturbation of observation matrix can be seen as an extension of the problem of solving systems of equations affected by uncertainty which has received significant attention in the literature (cf., e.g., \cite{cope1979bounds,higham2002accuracy,kreinovich1993optimal,nazin2005interval,neumaier1990interval,oettli1964compatibility,polyak2003robust} and references therein).
It is also tightly related to the system identification problem under uncertain-but-bounded perturbation of the observation of the state of the system \cite{bertsekas1971recursive,casini2014feasible,cerone1993feasible,JuKoNe,kurzhansky1997ellipsoidal,matasovestimators,milanese2013bounding,nazin2007ellipsoid,tempo1990optimal}.

In what follows, our goal is to extend the estimation constructions from \cite{PUP} to the case of uncertain sensing matrix.
Our strategy {consists in constructing} a tight efficiently computable convex in $H$ upper bound on the risk of
a candidate
 estimate, and then {building} a ``presumably
 good'' estimate by minimizing this bound in the estimate parameter $H$.
Throughout the paper, we assume that
 the signal set $\cX$ is an ellitope, and the norm $\|\cdot\|$ quantifying the recovery error is the maximum of a finite collection of Euclidean norms.
 \paragraph{Our contributions} can be summarized as follows.
\begin{itemize}
\item[A.] In Section \ref{secplgle} we analyse the $\epsilon$-risk
(the maximum, over signals from $\cX$, of the radii of $(1-\epsilon)$-confidence $\|\cdot\|$-balls) and the design of presumably good, in terms of this risk,
 linear estimates in the case of random uncertainty.
  \item[B.] In Section \ref{sec:unbblin},
  we build presumably good linear estimates in the case of {\em structured norm-bounded uncertainty} (cf. \cite[Chapter 7]{RO} and references therein),  thus extending the corresponding results of \cite{JuKoNe}.
 \end{itemize}
Developments in  A and B lead to novel computationally efficient techniques for designing presumably good
linear estimates for both random and uncertain-but-bounded perturbations.
\par
Analysis and design of polyhedral estimates under uncertainty in sensing matrix form the subject of Sections \ref{pgle} (random perturbations) and
\ref{sec:polyunbb} (uncertain-but-bounded perturbations). The situation here is as follows:
\begin{itemize}
\item[C.] The random perturbation case of the
{\em Analysis problem}
\begin{quote} {\sl Given contrast matrix $H$, find a
provably tight efficiently computable upper bound on $\epsilon$-risk of the associated estimate}
\end{quote}
is the subject of Section \ref{pgle}, where it is solved ``in the full range'' of our assumptions (ellitopic $\cX$, sub-Gaussian zero mean
$\eta$  and $\xi$). In contrast, the random perturbation case of the {\sl Synthesis problem} in which we want to minimize  the above bound w.r.t. $H$
turns out to be more involving---the bound to be optimized happens to be nonconvex in $H$. When there is no uncertainty in sensing matrix, this
difficulty can be somehow circumvented  \cite[Section 5.1]{PUP}; however, when uncertainty in sensing matrix is present, the strategy developed in \cite[Section 5.1]{PUP}
happens to work only when $\cX$ is an ellipsoid rather than a general-type ellitope. The corresponding developments are the subject
of Sections \ref{sectstra}, \ref{sec:end_stoch}, and \ref{sectimpl}.
\item[D.] In our context, analysis and design of polyhedral estimates under uncertain-but-bounded perturbations in the sensing matrix
appears to be the
most difficult; our very limited results on this subject form the subject of Section \ref{sec:polyunbb},
\end{itemize}

\paragraph{Notation and assumptions.}
We denote with $\|\cdot\|$ the norm on $\bR^\nu$ used to measure the estimation error. In what follows, $\|\cdot\|$ is a maximum of Euclidean norms
 \beq\label{eq0}
 \|u\|=\max_{\ell\leq L}\sqrt{u^TR_\ell u}
 \eeq
where $R_\ell \in \bS^\nu_+$, $\ell=1,...,L$, are given matrices with  ${\sum}_\ell R_\ell \succ 0$.

Throughout the paper, unless otherwise is explicitly stated, we assume that observation noise $\xi$ is zero-mean sub-Gaussian, $\xi\sim \SG(0,\sigma^2I)$, i.e.,
for all $t\in \bR^m$,
\be
\bE\left\{e^{t^T\xi}\right\}\leq \exp\left(\tfrac{\sigma^2}{2}{\|t\|_2^2}\right).
\ee{xi-sub}

\section{Random perturbations}
In this section we assume that uncertainty $\eta$ is sub-Gaussian, with parameters $0,I$, i.e.,
\be
\bE\left\{e^{t^T\eta}\right\}\leq \exp\left(\half {\|t\|_2^2}\right)\quad \forall t\in \bR^q.
\ee{eta-sub}
In this situation, given $\epsilon\in(0,1)$, 
we quantify the quality of recovery $\wh w(\cdot)$ of $w=Bx$ by its maximal over $x\in \cX$  {\em $\epsilon$-risk}
\be
\Risk_{\epsilon}
[\wh{w}|\cX]:=\sup_{x\in \cX}\inf\left\{\rho:\,\Prob_{\xi,\eta}\{\|Bx-\widehat{w}(A[\eta]x+\xi)\|>\rho\}\leq \epsilon\right\}
\ee{erisk}
(the radius of the smallest $\|\cdot\|$-ball centered at $\wh w(\omega)$ which covers $x$, uniformly over $x\in \cX$).

 \subsection{Design of presumably good linear estimate}\label{secplgle}
 \subsubsection{Preliminaries: ellitopes}\label{sec:pelli}
Throughout this section, we assume that the signal set $\cX$ is {\em a basic ellitope}.
Recall that, by definition \cite{JudNem2018,PUP}, a basic ellitope in $\bR^n$ is a set of the form
\be   \cX=\{x\in\bR^n:\,\exists t\in\cT:\,z^TT_k z\leq t_k,\,k\leq K\},
\ee{cX}
where $T_k\in\bS^n_+$, $T_k\succeq0$, ${\sum}_k T_k \succ0$, and $\cT\subset\bR^{K}_+$ is a convex compact set with a nonempty interior which is monotone: whenever $0\leq t'\leq t\in\cT$ one has $t'\in\cT$. We refer to $K$ as {\em ellitopic dimension} of $\cX$.
\par
Clearly, every basic ellitope is a convex compact set with nonempty interior which is symmetric w.r.t. the origin.
For instance, 
\\
{\bf A.} Bounded intersection $\cX$ of $K$ centered at the origin ellipsoids/elli\-p\-tic cylinders $\{x\in \bR^n:x^T T_kx\leq1\}$ [$T_k\succeq0$] is a {basic ellitope}:
\[
\cX=\{x\in \bR^n:\exists t\in\cT:=[0,1]^K: x^T T_kx\leq t_k,\,k\leq K\}
\]
In particular, the unit box $\{x\in \bR^n:\|x\|_\infty\leq1\}$ is a basic ellitope.\\
{\bf B.} A $\|\cdot\|_p$-ball in $\bR^n$ with $p\in[2,\infty]$ is a basic ellitope:
\[
\{x\in\bR^n:\|x\|_p\leq1\} =\big\{x:\exists t\in\cT=\{t\in\bR^n_+,\|t\|_{p/2}\leq 1\}:\underbrace{ x_k^2}_{x^T T_k x}\leq t_k,\,k\leq K\big\}.
\]
In the present context, our interest for ellitopes is motivated by their special relationship with the optimization problem
 \begin{equation}\label{quadprob}
\Opt_*(C)=\max\limits_{x\in\cX}x^T Cx,\,\,C\in\bS^n
\end{equation}
of maximizing a homogeneous quadratic form  over $\cX$.
As it is shown in \cite{PUP}, when $\cX$ is an ellitope, \rf{quadprob} admits ``reasonably tight'' efficiently computable upper bound. Specifically,
\begin{theorem}\label{2020Prop4.6} {\rm \cite[Proposition 4.6]{PUP}} Given ellitope {\rm \rf{cX}} and matrix $C$, consider the quadratic maximization problem {\rm (\ref{quadprob})}
along with its relaxation\footnote{Here and below, we use notation $\phi_\cS(\cdot)$ for the support function of a convex set $\cS\subset\bR^n$: for $y\in \bR^n$,
 $$
 \phi_\cS(y)=\sup_{u\in\cS} y^Ts.
 $$}
\begin{equation}\label{2020eq10}
\Opt(C)=\min_\lambda\left\{\phi_{\cT}(\lambda): \lambda\geq0,{{\sum}}_k\lambda_kT_k-C\succeq 0\right\}
\end{equation}
  The problem  is computationally tractable and solvable, and $\Opt(C)$ is an efficiently computable upper bound on $\Opt_*(C)$.
This upper bound is tight:
\[
\Opt_*(C)\leq \Opt(C)\leq3\ln(\sqrt{3}K)\Opt_*(C).
\]
\end{theorem}
\subsubsection{Tight upper bounding the risk of linear estimate}
Consider a linear estimate
 \[
 \widehat{w}^H(\omega)=H^T\omega\quad[H\in\bR^{m\times\nu}]
 \]
 \begin{proposition}\label{pr:lin_stoch}
In the setting of this section, synthesis of a presumably good linear estimate reduces to solving  the convex optimization problem
\be
\min\limits_{H\in\bR^{m\times\nu}}\mR[H]
\ee{linrand}
where
\begin{equation}\label{prop31}
\begin{array}{rcl}
\mR[H]&=&\min_{\lambda_\ell,\mu^\ell, \kappa^\ell,\atop\varkappa^\ell,\rho,\varrho}
\bigg\{\left[1+\sqrt{2\ln(2L/\epsilon)}\right]\left[\sigma\max\limits_{\ell\leq L}\|HR_\ell^{1/2}\|_\Fro+\rho\right]+\varrho:\\
&& \left.\begin{array}{l}\mu^\ell\geq0,\varkappa^\ell\geq0,\,
\lambda_\ell+\phi_\cT(\mu_\ell)\leq\rho,\kappa_\ell+\phi_\cT(\varkappa^\ell)\leq\varrho,\,\ell\leq L\\
\hbox{\scriptsize$\left[\begin{array}{c|c}\lambda_\ell I_{\nu q}&{1\over 2}\left[R_\ell^{1/2}H^TA_1;...;R_\ell^{1/2}H^TA_q\right]\cr\hline {1\over 2}\left[A_1^THR_\ell^{1/2},...,A_q^THR_\ell^{1/2}\right]&{\sum}_k\mu^\ell_kT_k\cr\end{array}\right]$}\succeq0,\,\ell\leq L\\
\hbox{\scriptsize$\left[\begin{array}{c|c}\kappa^\ell I_\nu&{1\over 2}R_\ell^{1/2}[B-H^TA]\cr\hline {1\over 2}[B-H^TA]^TR_\ell^{1/2}&{\sum}_k\varkappa^\ell_kT_k\cr\end{array}\right]$}\succeq0,\,\ell\leq L\\
\end{array}\right\}\\
\end{array}
\end{equation}
For a candidate contrast matrix $H$, the $\epsilon$-risk of the linear estimate $\widehat{w}^H_\lin(\omega)=H^T\omega$ is upper-bounded by $\mR[H]$.
\par
 \end{proposition}
\subsubsection{A modification}\label{sectmod}
Let us assume  that a $K$-repeated  version of observation  (\ref{eq1}) is available, i.e., we {observe}
\begin{equation}\label{eq1rep}
\omega^K=\{\omega_k=A[\eta_k]x+ \xi_k,\,k=1,...,K\}
\end{equation}
with independent across $k$ pairs $(\xi_k,\eta_k)$. In this situation, we can relax the assumption of sub-Gaussianity of $\xi$ and $\eta$ to the second moment boundedness condition
\begin{equation}\label{relax}
\bE\{\xi\xi^T\}\preceq \sigma^2 I_m,\quad\bE\left\{\eta\eta^T\right\}\preceq I_q.
\end{equation}
Let us consider the following construction.
For each $\ell\leq L$, given  $H\in \bR^{m\times\nu}$ we denote
\begin{equation}\label{add0}
\begin{array}{rcl}
\widetilde\mR_\ell[H]&=&\min_{\lambda,\mu, \kappa,\varkappa}
\bigg\{\sigma\|HR_\ell^{1/2}\|_\Fro+\lambda+\phi_\cT(\mu)+\kappa+\phi_\cT(\varkappa):\\
&& \left.\begin{array}{l}\mu\geq0,\varkappa\geq0,\left[\begin{array}{c|c}\kappa I_\nu&{1\over 2}R_\ell^{1/2}[B-H^TA]\cr\hline {1\over 2}[B-H^TA]^TR_\ell^{1/2}&{\sum}_k\varkappa_kT_k\cr\end{array}\right]\succeq0
\\
\left[\begin{array}{c|c}\lambda I_{\nu q}&{1\over 2}\left[R_\ell^{1/2}H^TA_1;...;R_\ell^{1/2}H^TA_q\right]\cr\hline {1\over 2}\left[A_1^THR_\ell^{1/2},...,A_q^THR_\ell^{1/2}\right]&{\sum}_k\mu_kT_k\cr\end{array}\right]\succeq0
\end{array}\right\}
\end{array}
\end{equation}
and consider the convex optimization problem
\begin{equation}\label{linexp}
\wt H_\ell\in \Argmin_H\widetilde\mR_\ell[H].
\end{equation}

We define the ``reliable estimate'' $\wh w^{(r)}(\omega^K)$ of $w=Bx$ as follows.
\begin{enumerate}\item Given $H_\ell\in \bR^{m\times\nu} $ and observations $\omega_k$ we compute linear estimates $w_{\ell}(\omega_k)=H_\ell\omega_k$, $\ell=1,...,L,\,k=1,...,K$;
\item We define vectors $z_\ell\in \bR^\nu$ as geometric medians of $ w_{\ell}(\omega_k)$:
\[
z_\ell(\omega^K)\in \Argmin_z\sum_{k=1}^K\|R_\ell^{1/2}(w_\ell(\omega_k)-z)\|_2,\;\ell=1,...,L.
\]
\item Finally, we select as $\wh w^{(r)}(\omega^K)$ any point of the set
\[
\cW(\omega^K)=\bigcap_{\ell=1}^L\left\{w\in \bR^{\nu}:\,\|R_\ell^{1/2}(z_\ell(\omega^K)-w)\|_2\leq 4\wt \mR_\ell[H_\ell]\right\}.
\]
or set $\wh w^{(r)}(\omega^K)$ a once for ever fixed point, e.g., $\wh w^{(r)}(\omega^K)=0$ if $\cW(\omega^K)=\emptyset$.
\end{enumerate}
We have the following analog of Proposition \ref{pr:lin_stoch}.
\begin{proposition}\label{proplin2} In the situation of this section, it holds
\be\sup_{x\in \cX}\bE_{\eta_k,\xi_k}\left\{\|R_\ell^{1/2}(w_{\ell}(\omega_k)-Bx)\|_2^2\right\}\leq \widetilde\mR^2_\ell[H_\ell],\;\ell\leq L,
\ee{expkl}
and
\be
\Prob\left\{\|R_\ell^{1/2}(z_\ell(\omega^K)-Bx)\|_2\geq 4\wt \mR_\ell[H_\ell]\right\}\leq e^{-0.1070 K},\;\ell\leq L.
\ee{prbl}
As a consequence, whenever $K\geq \ln[L/\epsilon]/0.1070$, the $\epsilon$-risk of the aggregated estimate $\wh w^{(r)}(\omega^K)$ satisfies
\[\Risk_{\epsilon}
[\wh w^{(r)}(\omega^K)|\cX]\leq \ov \mR,\;\;\ov \mR=8\max_{\ell\leq L}\wt \mR_\ell[H_\ell].
\]
\end{proposition}
\paragraph{Remark.} Proposition \ref{proplin2} is motivated by the desire to capture situations in which sub-Gaussian assumption on $\eta$ and $\xi$ does not hold or is too restrictive. Consider, e.g., the case where the uncertainty in the sensing matrix reduces to zeroing out some randomly selected columns in the nominal matrix $\overline{A}$ (think of taking picture through the window with frost patterns). Denoting by $\gamma$ the probability to zero out a particular column and assuming that columns are zeroed out independently, model (\ref{eq1}) in this situation reads
$$
\omega=A[\eta]x+\xi,\, A[\eta]=(1-\gamma)\overline{A}+{\sum}_{\alpha=1}^n\eta_\alpha A_\alpha
$$
where  $\eta_1,...,\eta_n$  are i.i.d. zero mean random variables taking values $(\gamma-1)\rho$ and $\gamma \rho$ with probabilities $\gamma$ and $1-\gamma$, and $A_\alpha$, $1\leq\alpha\leq n$, is an $m\times n$ matrix with all but the $\alpha$-th column being zero and $\Col_\alpha[A_\alpha]=\rho^{-1}\Col_\alpha[\overline{A}]$. Scaling factor $\rho$ is selected to yield {the unit sub-Gaussianity parameter of $\eta$} or $\bE\{\eta_\alpha^2\}=1$ depending on whether Proposition \ref{pr:lin_stoch} or Proposition \ref{proplin2}
is used. For small $\gamma$, the scaling factor $\rho$  is essentially smaller in the first case, resulting in larger ``disturbance matrices'' $A_\alpha$ and therefore---in {stricter} constraints in the optimization problem \rf{linrand}, \rf{prop31} responsible for the design of the linear estimate.
\subsubsection{Numerical illustration}
\begin{figure}[htb!]
\begin{center}{\small
\begin{tabular}{c}
\includegraphics[width=0.9\textwidth]{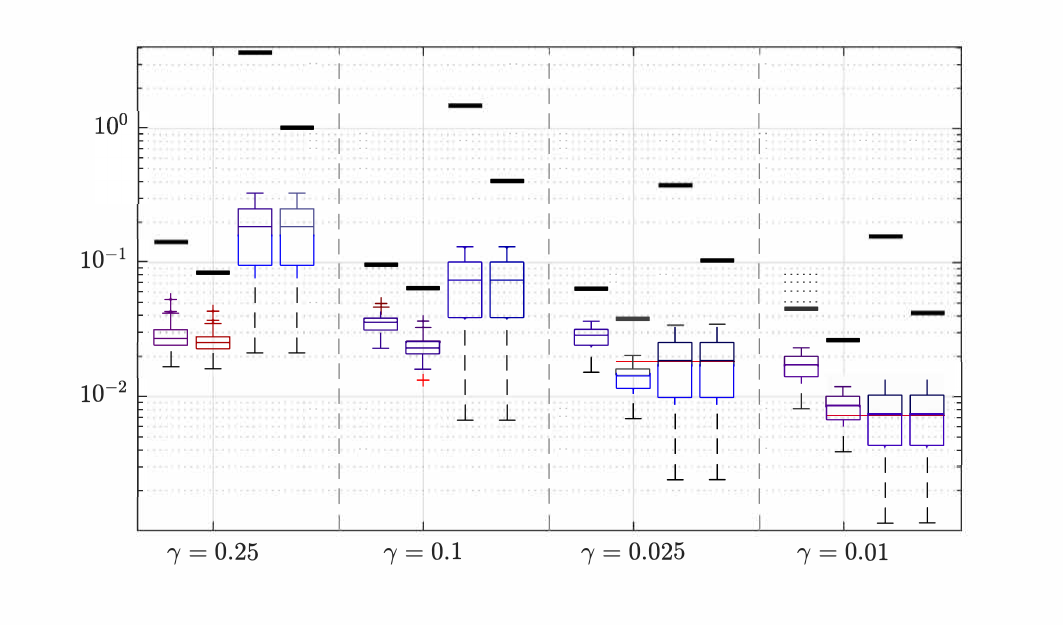}
\end{tabular}
\caption{\label{fig:1}  Distributions of $\ell_2$-recovery errors and upper bounds of the robust and ``nominal'' estimates for different values of $\gamma$ parameter.}
}\end{center}
\end{figure}
In Figure \ref{fig:1} we present results of a toy experiment in which
\begin{itemize}
\item $n=32,m=32$, and $\nu=16$,  $\overline{A}x\in\bR^m$ is the discrete time convolution of $x\in\bR^n$ with a simple kernel $\varkappa$ of length $9$ restricted onto the ``time horizon'' $\{1,...,n\}$, and $Bx$ cuts off $x$ the first $\nu$ entries. We consider Gaussian perturbation $\eta\sim\cN(0,\gamma^2 I_q)$,  $q=9$, and $A[\eta]x=[A+\sum_{\alpha=1}^q \eta_\alpha A_\alpha]x$ which is the convolution of $x$ with the kernel $\varkappa_\eta$ restricted onto the time horizon $\{1,...,n\}$, $\gamma$ being the control parameter.
\item $L=1$ and $\|\cdot\|=\|\cdot\|_2$,
\item $\cX$ is the ellipsoid $\{x: \sum_ii^2[Dx]_i^2\leq1\}$, where $D$ is the matrix of
    inverse Discrete Cosine Transform of size $n\times n$.
\item $\xi\sim\cN(0,\sigma^2 I_m)$, $\sigma=10^{-4}$.
\end{itemize}
In each cell of the plot we represent error distributions and upper risk bounds (horizontal bar) of four estimates (from left to right) for different uncertainty levels $\gamma$: (1) robust estimate by Proposition \ref{pr:lin_stoch} and upper bound $\mR$ on its $0.05$-risk, (2) single-observation estimate $w_1(\omega_1)=H_1\omega_1$  yielded by the minimizer $H_1$ of $\wt\mR_1[H]$  over $H$, see (\ref{add0}), and upper bound $\wt\mR_1[H_1]$ on its {\em expected error risk},\footnote{We define expected error risk of a $K$-observation estimate $\widehat{x}(\omega^K)$ of $Bx$ as $\sup_{x\in\cX}\bE_{\omega^K\sim P^K_x}\{\|\widehat{x}(\omega^K)-Bx\|\}$, where $P^K_x$ is the distribution of $\omega^K$  stemming from $x$.} (3) ``nominal'' estimate---estimate by
Proposition \ref{pr:lin_stoch} as applied to the ``no uncertainty'' case where all $A_\alpha$ in (\ref{eq2}) are set to 0 and upper bound $\mR$ from \rf{prop31} on its $0.05$-risk computed using actual uncertainty level,
(4) ``nominal'' estimate $\wt w_1(\omega_1)=\wt H_1\omega_1$ yielded by the minimizer $\wt H_1$ of $\wt\mR_1[H]$  over $H$ in the ``no uncertainty'' case and upper bound {$\wt \mR_1[\wt H_1]$}  on its ``actual''---with uncertainty present---expected error risk.

\subsection{Design of presumably good polyhedral estimate}\label{pgle}
\subsubsection{Preliminaries on polyhedral estimates} Consider a slightly more general than (\ref{eq1}), (\ref{eq2}) observation scheme
\beq\label{obseq}
\omega=Ax+\zeta
\eeq
where $A\in\bR^{m\times n}$ is given, unknown signal $x$ is known to belong to a given signal set $\cX$ given by (\ref{cX}), and $\zeta$ is observation noise with probability distribution $P_x$ which can depend on $x$.
For example, when observation $\omega$ is given by (\ref{eq1}), (\ref{eq2}), we have
\beq\label{wehave}
\zeta={\sum}_{\alpha=1}^q\eta_\alpha A_\alpha x +\xi
\eeq with  zero mean sub-Gaussian $\eta$ and $\xi$. \par
When building polyhedral estimate in the situation in question, one, given tolerance $\epsilon\in(0,1)$ and a positive integer $M$, specifies a computationally tractable convex set $\cH$, the larger the better, of vectors $h\in\bR^m$ such that
\beq\label{eq10}
\Prob_{\zeta\sim P_x}\{|h^T\zeta|>1\}\leq \epsilon/M\quad \forall x\in\cX.
\eeq
A polyhedral estimate $\widehat{w}^H(\cdot)$ is specified by {\sl contrast matrix} $H\in\bR^{M\times n}$ restricted to have all columns in $\cH$ according to
\beq\label{Au}
\omega\mapsto\wh x^H(\omega)\in\Argmin_{u\in\cX}\left\{\|H^T[Au-\omega]\|_\infty\right\},\;\widehat{w}^H_\poly(\omega)=B\wh x^H(\omega).
\eeq
It is easily seen (cf. \cite[Proposition 5.1.1]{PUP}) that the $\epsilon$-risk \rf{erisk}
of the above estimate is upper-bounded by the quantity
\begin{equation}\label{easily}
\mP[H]=\sup_y\left\{\|By\|:y\in 2\cX, \|H^TAy\|_\infty\leq 2\right\}.\
\end{equation}
\begin{quote}
Indeed, let $h_1,...,h_M$ be the columns of $H$. For $x\in\cX$ fixed, the inclusions  $h_j\in\cH$ imply that the $P_x$-probability of the event $Z_x=\{\zeta:|\zeta^Th_j|\leq1\,\forall j\leq M\}$ is at least $1-\epsilon$. When this event takes place, we have $\|H^T[\omega -Ax]\|_\infty\leq 1$, which combines with $x\in\cX$ to imply that $\|H^T[\omega-A\wh x^H(\omega)]\|_\infty\leq1$, so that $\|H^TA[x-\wh x^H(\omega)]\|_\infty\leq 2$, and besides this, $x-\wh x^H(\omega)\in 2\cX$, whence $\|Bx-
\widehat{w}^H_\poly(\omega)\|\leq \mP[H]$ by definition of $\mP[H]$. The bottom line is that whenever $x\in\cX$ and $\zeta=\omega-Ax\in Z_x$, which happens with $P_x$-probability at least $1-\epsilon$, we have $\|Bx-
\widehat{w}^H_\poly(\omega)\|\leq \mP[H]$, whence the $\epsilon$-risk of the estimate $\widehat{w}^H_\poly$ indeed is upper-bounded by $\mP[H]$.
\end{quote}
To get a presumably good polyhedral estimate, one minimizes $\mP[H]$ over $M\times \nu$ matrices $H$ with columns from $\cH$. Precise minimization is problematic, because $\mP[\cdot]$, while being convex, is usually difficult to compute. Thus, the design routine proposed in \cite{juditsky2020polyhedral} goes via minimizing an efficiently computable upper bound on $\mP[H]$.
It is shown in \cite[Section 5.1.5]{PUP} that when $\cX$ is ellitope (\ref{cX}) and $\|u\|=\| Ru\|_2$, a reasonably tight upper bound on $\mP[H]$ is given by the efficiently computable function
$$
\mP_+[H]=2\min\limits_{\lambda,\mu,\upsilon}\left\{\lambda+\phi_\cT(\mu)+{\sum}_i\upsilon_i:\begin{array}{l}\mu\geq0,\upsilon\geq0\\
\hbox{\scriptsize$\left[\begin{array}{c|c}\lambda I_\nu&{1\over 2}RB\cr\hline{1\over 2}B^TR^T&A^TH\Diag\{\upsilon\}H^TA+{\sum}_k\mu_kT_k\cr\end{array}\right]\succeq0$}
\end{array}\right\}.
$$
Synthesis of a presumably good polyhedral estimate reduces to minimizing the latter function in $H$ under the restriction $\Col_j[H]\in\cH$. Note that the latter problem still is nontrivial because $\mP_+$ is nonconvex in $H$.
\par
Our objective here is to implement the outlined strategy in the case of observation $\omega$ is given by \rf{eq1}, \rf{eq2}.
\subsubsection{Specifying $\cH$}\label{sectcH}
Our first goal is to specify, given tolerance $\delta\in(0,1)$, a set $\cH_\delta\subset\bR^m$, the larger the better, such that
\beq\label{goal1}
h\in\cH_\delta, x\in\cX \Rightarrow \Prob_{\zeta\sim P_x}\{|h^T\zeta|>1\} \leq \delta.
\eeq
Note that a ``tight'' sufficient condition for the validity of (\ref{goal1}) is
\begin{subequations}
\label{condab}
\begin{align}\label{condab.a}
\Prob_{\xi}\{|h^T\xi|>1/2\}&\leq\delta/2,\\
\label{condab.b}
\Prob_{\eta}\left\{\left|{\sum}_{\alpha=1}^q [h^TA_\alpha x]\eta_\alpha\right|>1/2\right\}&\leq\delta/2,\,\forall x\in\cX.
\end{align}
\end{subequations}
Note that under the sub-Gaussian assumption \rf{xi-sub}, $h^T\xi$ is itself sub-Gaussian, $h^T\xi\sim \SG(0,\sigma^2\|h\|_2^2)$; thus, a tight sufficient condition for \rf{condab.a} is
\beq\label{conda}
\|h\|_2\leq [\sigma\chi(\delta)]^{-1},\,\,\chi(\delta)=2\sqrt{2\ln(2/\delta)}.
\eeq
Furthermore, by \rf{eta-sub}, r.v. ${\sum}_{\alpha=1}^q [h^TA_\alpha x]\eta_\alpha=h^T[A_1x,...,A_qx]\eta$ is sub-Gaussian with parameters $0$ and
$\| [h^TA_1x;...;h^TA_qx]\|_2^2$, implying the validity of \rf{condab.b} for a given $x$ whenever
\[
\|[h^TA_1x;...;h^TA_qx]\|_2\leq \chi^{-1}(\delta).
\]
We want this relation to hold true for every $x\in\cX$, that is, we want the operator norm $\|\cdot\|_{\cX,2}$ of the mapping
\beq\label{barA}
x\mapsto \cA[h]x,\;\cA[h]=[h^TA_1;h^TA_2;...;h^TA_q]
\eeq
induced by the norm $\|\cdot\|_\cX$ on the argument and the norm $\|\cdot\|_2$ on the image space to be upper-bounded by $\chi(\delta)$:
\be
\|\cA[h]\|_{\chi,2}\leq \chi^{-1}(\delta).
\ee{condb}
 Invoking \cite[Theorem 3.1]{JuKoNe} (cf. also the derivation in the proof of Proposition \ref{pr:lin_stoch} in Section \ref{sec:proof31}), a tight sufficient condition for the latter relation is
\beq\label{eq100}
\Opt[h]:=\min_{\lambda,\mu}\left\{\lambda+\phi_\cT(\mu):\,\mu\geq0,\,\left[\begin{array}{c|c}\lambda I_q&{\half}\cA[h]\cr\hline{\half}\cA^T[h]&{\sum}_k\mu_kT_k\cr\end{array}\right]\succ 0\right\}\leq \chi^{-1}(\delta),
\eeq
tightness meaning that $\Opt[h]$ is within factor $O(1)\sqrt{\ln(K+1)}$ of $\|\cA[h]\|_{\cX,2}$.
\par
The bottom line is that with $\cH_\delta$ specified by constraints (\ref{conda}) and \rf{condb} (or by the latter replaced with its tight relaxation \rf{eq100}) we do ensure (\ref{goal1}).
\subsubsection{Bounding the risk of the polyhedral estimate $\wh w^H$}\label{sect323}
\begin{proposition}\label{propanal}
In the situation of this section, let $\epsilon\in (0,1)$, and let $H=[H_1,...,H_L]$ be $m\times ML$ matrix with $L$ blocks $H_\ell\in\bR^{m\times M}$  such that $\Col_j[H]\in\cH_\delta$ for all $j\leq ML$ and $\delta=\epsilon/ML$.  Consider optimization problem
\begin{align}\label{optprob}
\mP_+[H]&=2\min_{\lambda_\ell,\mu^\ell,\upsilon^\ell,\rho}\left\{
\rho:\,\mu^\ell\geq0,\upsilon^\ell\geq0, \,\lambda_\ell+\phi_\cT(\mu^\ell)+{\sum}_{j=1}^M\upsilon^\ell_j \leq\rho,\,\ell\leq L\right.\nn
&\qquad\qquad\qquad\left.\left[\begin{array}{c|c}\lambda_\ell I_\nu&{1\over 2}R_\ell^{1/2}B\cr\hline {1\over 2}B^TR_\ell^{1/2}&A^TH_\ell\Diag\{\upsilon^\ell\}H_\ell^TA+
{\sum}_k\mu^\ell_kT_k\cr\end{array}\right]\succeq0,\,\ell\leq L\right\}.
\end{align}
Then
\[
\Risk_{\epsilon}[\widehat{w}^H|\cX]\leq \mP_+[H].
\]
\end{proposition}

\subsubsection{Optimizing $\mP_+[H]$---the strategy}\label{sectstra}
Proposition \ref{propanal} resolves the {\sl analysis} problem---it allows to efficiently upper-bound the $\epsilon$-risk of a given polyhedral estimate $\wh w^H_\poly$. At the same time, ``as is,'' it does not allow to build the estimate itself (solve the ``estimate synthesis'' problem---compute a presumably good contrast matrix) because straightforward minimization of $\mP_+[H]$ (that is,
adding $H$ to decision variables of the right hand side of (\ref{optprob}) results in a nonconvex problem. A remedy,
as proposed in \cite[Section 5.1]{PUP}, stems from the concept  of a {\sl cone compatible with a convex compact set $\cH\subset\bR^m$} which is defined as follows:
\par
Given positive integer $J$ and real $\varkappa\geq1$ we say that a closed convex cone $\bK\subset\bS^m_+\times\bR_+$ is $(J,\varkappa)$-compatible with $\cH$
if
\begin{enumerate}
\item[(i)] whenever $h_1,...,h_J\in\cH$ and $\upsilon\in\bR^J_+$, the pair $\left({\sum}_{j=1}^J\upsilon_j h_jh_j^T,{\sum}_j\upsilon_j\right)$ belongs to $\bK$,
\\
 and ``nearly vice versa'':
\item[(ii)] given $(\Theta,\varrho)\in\bK$ {and $\varkappa\geq 1$}, we can efficiently build collections of vectors $h_j\in\cH$, and reals $\upsilon_j\geq0$,
$j\leq J$, such that $\Theta={\sum}_{j=1}^J\upsilon_jh_jh_j^T$ and ${\sum}_j\upsilon_j\leq \varkappa \varrho$.
\end{enumerate}
{\bf Example.} Let $\cH$ be a centered at the origin Euclidean ball of radius $R>0$ in $\bR^J$. When setting \[
\bK=\{(\Theta,\varrho):\,\Theta\succeq0,\Tr(\Theta)\leq R^2\varrho\Tr(\Theta)\},
\]
 we obtain a cone $(M,1)$-compatible with $\cH$. Indeed, for
$h_j\in\cH$ and $\upsilon_j\geq0$ we have
\[\Tr\left({\sum}_j\upsilon_jh_jh_j^T\right)\leq R^2{\sum}_j\upsilon_j,
\] that is $\left(\Theta:={\sum}_j\upsilon_jh_jh_j^T,\varrho:={\sum}_j{\upsilon}_j\right)\in\bK$. Vice versa, given $(\Theta,\varrho)\in\bK$, i.e., $\Theta\succeq0$ and $\varrho\geq\Tr(\Theta)/R^2$ and specifying $f_1,...,f_m$ as the orthonormal system of eigenvectors of $\Theta$, and $\lambda_j$ as the corresponding eigenvalues and setting $h_j=Rf_j$, $\upsilon_j=R^{-2}\lambda_j)$, we get $h_j\in\cH$,
$\Theta={\sum}_j\upsilon_jh_jh_j^T$ and ${\sum}_j\upsilon_j=\Tr(\Theta)/R^2\leq\varrho$.
\par
Coming back to the problem of minimizing $\mP_+[H]$ in $H$,  assume that we have at our disposal a cone $\bK$ which is $(M,\varkappa)$-compatible with $\cH_\delta$. In this situation, we can replace the nonconvex problem
\be\min_{H=[H^1,...,H^L]} \{\mP_+[H]:\,\Col_j[H^\ell]_j\in \cH_\delta\}
\ee{proboptH}
with the problem
\begin{align}\label{proboptHH}
&\min\limits_{{\bar{\lambda}_\ell,\bar{\mu}^\ell,\atop \Theta_\ell,\varrho_\ell,\bar{\rho}}}
\Big\{
\bar{\rho}:\,(\Theta_\ell,\varrho_\ell)\in\bK,\bar{\mu}^\ell\geq0,\,\bar{\lambda}_\ell+\phi_\cT(\bar{\mu}^\ell)+\varrho_\ell \leq\bar{\rho},
\,\ell\leq L,\nn
&\qquad\qquad\qquad\left.\left[\begin{array}{c|c}\bar{\lambda}_\ell I_\nu&{1\over 2}R_\ell^{1/2}B\cr\hline {1\over 2}B^TR_\ell^{1/2}&A^T\Theta_\ell A+
{\sum}_k\bar{\mu}^\ell_kT_k\cr\end{array}\right]\succeq0, \ell\leq L\right\}.
\end{align}
Unlike \rf{proboptH}, the latter problem is convex and efficiently solvable provided that $\bK$ is computationally tractable, and can be considered as ``tractable $\sqrt{\varkappa}$-tight'' relaxation of the problem of interest (\ref{proboptH}). Namely,
\begin{itemize}
\item Given a feasible solution $H_\ell,\lambda_\ell,\mu^\ell,\upsilon^\ell,\rho$ to the problem of interest (\ref{proboptH}), we can set \[\Theta_\ell={\sum}_{j=1}^M\upsilon^\ell_j\Col_j[H_\ell]\Col^T_j[H_\ell], \;\;{\varrho}_\ell={\sum}_j\upsilon^\ell_j,
     \]thus getting $(\Theta_\ell,{\varrho}_\ell)\in\bK$. By (i) in the definition of compatibility, $\Theta_\ell,\varrho_\ell,\bar{\lambda}_\ell=\lambda_\ell,\bar{\mu}^\ell=\mu^\ell,\bar{\rho}=\rho$ is a feasible solution to (\ref{proboptHH}), and this transformation preserves the value of the objective
\item Vice versa, given a feasible solution $\Theta_\ell,\varrho_\ell,\bar{\lambda}_\ell,\bar{\mu}^\ell,\bar{\rho}$ to (\ref{proboptHH}) and invoking (ii) of the definition of compatibility, we can convert, in a computationally efficient way, the pairs $(\Theta_\ell,\rho_\ell)\in\bK$ into the pairs $H_\ell\in\bR^{m\times M}$, $\bar{\upsilon}^\ell\in\bR^m_+$ in such a way that the columns of $H_\ell$ belong to $\cH_\delta$, $\Theta_\ell=H_\ell\Diag\{\bar{\upsilon}^\ell\}H_\ell^T$, ${\sum}_j\bar{\upsilon}^\ell_j\leq\varkappa\varrho_\ell$. Assuming w.l.o.g. that all matrices $R_\ell^{1/2}B$ are nonzero, we obtain $\phi_\cT(\bar{\mu}^\ell)+\varrho_\ell>0$ and $\bar{\lambda}_\ell>0$ for all $\ell$. We claim that setting
$$
\gamma_\ell=\sqrt{[\phi_\cT(\bar{\mu}^\ell)+\varkappa\varrho_\ell]/\bar{\lambda}_\ell},\;
\lambda_\ell=\gamma_\ell\bar{\lambda}_\ell,\;\mu_\ell=\gamma_\ell^{-1}\bar{\mu}_\ell,\upsilon^\ell=\gamma_\ell^{-1}\bar{\upsilon}^\ell,\;
\rho=\sqrt{\varkappa}\bar{\rho}
$$
we get a feasible solution to (\ref{proboptH}). Indeed, all we need is to verify that this solution satisfies, for every $\ell\leq L$,  constraints  of \rf{optprob}.
To check the semidefinite constraint, note that
$$
\hbox{\footnotesize$ \left[\begin{array}{c|c}\lambda_\ell I_\nu&{1\over 2}R_\ell^{1/2}B\cr\hline
{1\over 2}B^TR_\ell^{1/2}&A^TH_\ell\Diag\{\upsilon^\ell\}H_\ell^TA+{\sum}_k\mu^\ell_kT_k\cr\end{array}\right]
=\left[\begin{array}{c|c}\gamma_\ell\bar{\lambda}_\ell I_\nu&{1\over 2}R_\ell^{1/2}B\cr\hline
{1\over 2}B^TR_\ell^{1/2}&\gamma_\ell^{-1}\left[A^TH_\ell\Diag\{\bar{\upsilon}^\ell\}H_\ell^TA+{\sum}_k\bar{\mu}^\ell_kT_k\right]\cr\end{array}\right]
$}$$
and the matrix in the right-hand side is $\succeq0$ by the semidefinite constraint of \rf{proboptHH} combined with $\Theta_\ell={\sum}_j\bar{\upsilon}^\ell_j\Col_j[H_\ell]\Col^T_j[H_\ell]$. Furthermore, note that by construction ${\sum}_j\bar{\upsilon}^\ell_j\leq\varkappa\varrho_\ell$, whence
\begin{align*}
\lambda_\ell+\phi_\cT(\mu^\ell)+{\sum}_j\upsilon^\ell_j
&=\gamma_\ell\bar{\lambda_\ell}+\gamma_\ell^{-1}[\phi_\cT(\bar{\mu}^\ell)+\varkappa\varrho_\ell]=2\sqrt{\bar{\lambda}_\ell
[\phi_\cT(\bar{\mu}^\ell)+\varkappa\varrho_\ell]}\\
&\leq2\sqrt{\varkappa}\sqrt{\bar{\lambda}_\ell
[\phi_\cT(\bar{\mu}^\ell)+\varrho_\ell]}\leq \sqrt{\varkappa}\left[\bar{\lambda}_\ell+\phi_\cT(\bar{\mu}^\ell)+\varrho_\ell\right]\leq\sqrt{\varkappa}\bar{\rho}=\rho
\end{align*}
(we have taken into account that $\varkappa\geq1$).
\end{itemize}
We conclude that
the (efficiently computable) optimal solution to the relaxed problem (\ref{proboptHH}) can be efficiently converted to
a feasible solution to problem (\ref{proboptH}) which is within the factor at most $\sqrt{\varkappa}$  from optimality in terms of the objective. Thus,
\begin{quote}
(!) {\sl Given a $\varkappa$-compatible with $\cH_\delta$ cone $\bK$, we can find, in a computationally efficient fashion, a feasible solution to the problem of interest (\ref{proboptH}) with the value of the objective by at most the factor $\sqrt{\varkappa}$ greater than the optimal value of the problem.}
\end{quote}
\par
What we propose is to build a presumably good polyhedral estimate by applying the just outlined strategy to the instance of (\ref{proboptH}) associated with $\cH=\cH_\delta$ given by (\ref{conda}) and {\rf{eq100}}. The still missing---and crucial---element in this strategy is
a computationally tractable cone $\bK$ which is $(M,\varkappa)$-compatible, for some ``moderate'' $\varkappa$, with our $\cH_\delta$. For the time being,  we have at our disposal such a cone only for the ``no uncertainty in sensing matrix'' case (that is, in the case where all $A_\alpha$ are zero matrices), and  it  is shown in \cite[Chapter 5]{PUP} that in this case the polyhedral estimate stemming from the just outlined strategy  is near minimax-optimal, provided that $\xi\sim\cN(0,\sigma^2I_m)$.
\par
When ``tight compatibility''---with $\varkappa$ logarithmic in the dimension of $\cH$---is sought, the task of building a cone $(M,\varkappa)$-compatible with a given convex compact set $\cH$ reveals to be highly nontrivial. To the best of our knowledge, for the time being, the widest family of sets $\cH$ for which tight compatibility has been achieved is the family of ellitopes \cite{InProgr}. Unfortunately, this family seems to be too narrow to capture the sets $\cH_\delta$ we are interested in now. At present, the only known to us ``tractable case''  here is the ball case $K=1$, and even handling this case requires extending compatibility results of \cite{InProgr} from ellitopes to {\sl spectratopes.}

\subsubsection{Estimate synthesis utilizing cones compatible with spectratopes}
\label{sec:end_stoch}
Let for $S^{ij}\in\bS^{d_i}$, $1\leq i\leq I$, $1\leq j\leq N$, and let for $g\in \bR^N$,  $S_i[g]={\sum}_{j=1}^Ng_jS^{ij}$. A {\em basic spectratope in $\bR^N$} is a set  $\cH\subset\bR^N$ represented as
\beq\label{def+spec}
\cH=\{g\in\bR^N: \exists r\in\cR: S_i^2[g]\preceq r_iI_{d_i},i\leq I\};
\eeq
here $\cR$ is a compact convex {\sl monotone} subset of $\bR^I_+$ with nonempty interior, and ${\sum}_iS_i^2[g]\succ0$ for all $g\neq0$. We refer to $d={\sum}_i d_i$ as {\em spectratopic dimension} of $\cH$.
A spectratope, by definition, is a linear image of a basic spectratope.

As shown in \cite{PUP}, where the notion of a spectratope was introduced, spectratopes are convex compact sets symmetric w.r.t. the origin, and basic spectratopes have nonempty interiors. The family of spectratopes is rather rich---finite intersections, direct products, linear images, and arithmetic sums of spectratopes, same as inverse images of spectratopes under linear embeddings, are spectratopes, with spectratopic representations of the results readily given by spectratopic representations of the operands. \par
Every ellitope is a spectratope. An example of spectratope which is important to us is the set $\cH_\delta$ given by (\ref{conda}) and (\ref{condb}) in the ``ball case'' where $\cX$ is an ellipsoid (case of $K=1$). In this case, by one-to-one linear parameterization of signals
$x$, accompanied for the corresponding updates in $A,A_\alpha$, and $B$, we can assume that $T_1=I_n$ in \rf{cX}, so that $\cX$ is the unit Euclidean ball,
$$
\cX=\{x\in\bR^n: x^Tx\leq1\}.
$$
In this situation, denoting by $\|\cdot\|_{2,2}$ the spectral norm of a matrix, constraints (\ref{conda}) and (\ref{condb}) specify the set

\begin{equation}\label{eeeee}
\begin{array}{rcl}
\cH_\delta&=&\Big\{h\in\bR^m: \|h\|_2\leq{(\sigma\chi(\delta))^{-1}}, \|\cA[h]\|_{2,2}\leq{\chi^{-1}(\delta)}\Big\}\\
&=&\Big\{h\in\bR^m:\;\exists r\in \cR: S_j^2[h]\preceq r_jI_{d_j},j\leq 2\Big\}
\end{array}
\end{equation}
where $\cR=\{[r_1;r_2]:0\leq r_1,r_2\leq 1\}$,
\[
S_1[h]=\sigma\chi(\delta)\left[\begin{array}{c|c}&h\cr\hline h^T&\cr\end{array}\right]\in\bS^{m+1},
\;S_2[h]=\chi(\delta) \left[\begin{array}{c|c}&\cA[h]\cr\hline \cA[h]^T&\cr\end{array}\right]\in\bS^{m+q}
\]
with $d_1=m+1$, $d_2=m+q$.
We see that in the ball case $\cH_\delta$ is a basic spectratope.
\par
We associate with a spectratope $\cH$, as defined in \rf{def+spec}, linear mappings
$$
\cS_i[G]={\sum}_{p,q}G_{pq}S^{ip}S^{iq}:\bS^N\to \bS^{d_i}.
$$
Note that
\[
\cS_i\left[{\sum}_jg_jg_j^T\right]={\sum}_jS_i^2[g_j],\,\,g_j\in{\bR^N},
\]
and
\begin{subequations}\label{when1}
\begin{align}\label{when1a}
&G\preceq G'\Rightarrow \,\cS_i[G]\preceq \cS_i[G'],\\
\label{when1b}
&\{G\succeq0 \ \&\  \cS_i[G]=0\,\forall \ell\}\Rightarrow \,G=0.
\end{align}
\end{subequations}
A cone ``tightly compatible'' with a basic spectratope is given by the following
   \begin{proposition}\label{prop2} Let $\cH\subset\bR^N$ be a basic spectratope
\[
\cH=\{g\in\bR^N:\,\exists r\in\cR: \,S_i^2[g]\preceq r_i I_{d_i},\,i\leq I\}
\]
with ``spectratopic data'' $\cR$ and $S_i[\cdot]$, $i\leq I$, satisfying the requirements in the above definition.

Let us specify the closed convex cone $\bK\subset\bS^{N}_+\times\bR_+$ as
\[
\bK=\big\{(\Sigma,\rho)\in\bS^{N}_+\times\bR_+:\,\exists r\in\cR:\,
\cS_i[\Sigma]\preceq\rho r_iI_{d_i},i\leq I\big\}.
\]
Then
\item[\rm{(i)}] whenever $\Sigma={\sum}_j\lambda_jg_jg_j^T$ with $\lambda_j\geq0$ and $g_j\in\cH$ $\forall j$, we have
$$
\left(\Sigma,\,{{\sum}}_j\lambda_j\right)\in \bK,
$$
\item[\rm{(ii)}] and ``nearly'' vice versa: when $(\Sigma,\rho)\in\bK$, there exist (and can be found efficiently by a randomized algorithm)  $\lambda_j\geq0$ and $g_j$, $j\leq N$, such that
$$
\Sigma={\sum}_j\lambda_jg_jg_j^T\;\;\mathrm{with}\;\;{\sum}_j\lambda_j\leq\varkappa\rho\;\;\mathrm{and}\;\;g_j\in\cH,\;j\leq N.
$$
where
$$
\varkappa=4\ln(4DN),\,D={\sum}_i d_i.
$$
\end{proposition}
For the proof and for the sketch of the randomized algorithm mentioned in (ii), see Section \ref{sec:pr2pr} of the appendix.

\subsubsection{Implementing the strategy}\label{sectimpl}
We may now summarize our approach to the design of a presumably good polyhedral estimate. By reasons outlined at the end of Section \ref{sectstra}, the only case where the components we have developed so far admit ``smooth assembling'' is the one where $\cX$ is ellipsoid which in our context w.l.o.g. can be assumed to be  the unit Euclidean ball. Thus, {\sl in the rest of this Section it is assumed that $\cX$ is  the unit Euclidean ball in $\bR^n$.}
Under this assumption the recipe, suggested by the preceding analysis,  for designing presumably good polyhedral estimate is as follows. Given $\epsilon\in(0,1)$, we
\\$\bullet$ set $\delta=\epsilon/Lm$  and solve the convex optimization problem
\begin{align}
\label{final}
\Opt&=\min\limits_{{\Theta_\ell\in\bS^m,\atop
\varrho_\ell,\bar{\lambda}_\ell,\bar{\mu}_\ell}}\Big\{
\bar{\rho}:\,\bar{\mu}_\ell\geq0,\,\Theta_\ell\succeq0,\,\sigma^2\chi^2(\delta)\Tr(\Theta_\ell)\leq\varrho_\ell,\,
\bar{\lambda}_\ell+\bar{\mu}_\ell+\varrho_\ell \leq\bar{\rho},\,\ell\leq L,\nn
&\qquad\quad\left.\begin{array}{l}
\left[\begin{array}{c|c}\left[\Tr(A_\alpha^T\Theta_\ell A_\beta)\right]_{\alpha,\beta=1}^q&\cr\hline
&\sum\limits_{\alpha,\beta}A_\alpha^T\Theta_\ell A_\beta\cr\end{array}\right]\preceq\chi^{-2}(\delta)\varrho_\ell I_{q+n},\,\ell\leq L,\\
\left[\begin{array}{c|c}\bar{\lambda}_\ell I_\nu&
{1\over 2}R_\ell^{1/2}B\cr\hline {1\over 2}B^TR_\ell^{1/2}&A^T\Theta_\ell A+
\bar{\mu}_\ell I_n\cr\end{array}\right]\succeq0,\ell\leq L
\end{array}\right\}
\end{align}
---this is what under the circumstances becomes problem (\ref{proboptHH}) with the cone $\bK$ given by Proposition \ref{prop2} as applied to the spectratope $\cH_\delta$ given by (\ref{eeeee}). Note that by Proposition \ref{prop2}, $\bK$ is $\varkappa$-compatible with $\cH_\delta$, with
\begin{equation}\label{varkappa}
\varkappa=4\ln(4m(m+n+q+1)).
\end{equation}
For instance, in the case of rank 1 matrices
$A_\alpha=f_\alpha g_\alpha^T$ and $\|\cdot\|=\|\cdot\|_2$ \rf{final} becomes
\begin{align}
\label{finall}
\Opt&=\min_{\Theta\in\bS^m,\atop
\varrho,\bar{\lambda},\bar{\mu}}\Big\{
\bar{\rho}:\,\bar{\mu}\geq0,\,\Theta\succeq0,\,\sigma^2\chi^2(\delta)\Tr(\Theta)\leq\varrho,\,
\bar{\lambda}+\bar{\mu}+\varrho \leq\bar{\rho}\nn
&\qquad\quad\left.\begin{array}{l}
\left[\begin{array}{c|c}\left[(f_\alpha^T\Theta f_\beta)g_\alpha^Tg_\beta\right]_{\alpha,\beta=1}^q&\cr\hline
&\left[{\sum}_{\alpha,\beta=1}^q[f_\alpha^T\Theta f_\beta\right]g_\alpha g_\beta^T\cr\end{array}\right]\preceq\chi^{-2}(\delta)\varrho I_{q+n}\\
\left[\begin{array}{c|c}\bar{\lambda}_\ell I_\nu&
{1\over 2}B\cr\hline {1\over 2}B^T&A^T\Theta A+
\bar{\mu} I_n\cr\end{array}\right]\succeq0
\end{array}\right\};
\end{align}
$\bullet$ use the randomized algorithm described in the proof of Proposition \ref{prop2} to convert the $\Theta_\ell$-components of the optimal solution to (\ref{final}) into a contrast matrix. Specifically,
\begin{enumerate}\item
for $\ell=1,2,...,L$ we generate matrices $G_\varsigma^k=\Theta_\ell^{1/2}\Diag\{\varsigma^k\}O$, $k=1,...,K$, where $O$ is the orthonormal matrix of $m\times m$ Discrete Cosine Transform, and $\varsigma^k$ are i.i.d. realizations of $m$-dimensional Rademacher random vector;
\item for every $k\leq K$, we compute the maximum $\theta(G_\ell^k)$ of values of the Minkowski function of $\cH_\delta$ as evaluated
at the columns of $G_\ell^k$, with $\cH_\delta$  given by (\ref{conda}), (\ref{condb}), and select among $G_\ell^k$ matrix $G_\ell$ with the smallest value of $\theta(G_\ell^k)$.\\
Then the $\ell$-th block of the contrast matrix we are generating is $H_\ell=G_\ell\theta^{-1}(G_\ell)$.
\end{enumerate}
With reliability $1-2^{-K}L$ the resulting contrast matrix $H$ (which definitely has all columns in $\cH_\delta$) is, by (!), near-optimal, within factor $\sqrt{\varkappa}$ in terms of the objective, solution to (\ref{proboptH}), and the $\epsilon$-risk of the associated polyhedral estimate is upper-bounded by $2\sqrt{\varkappa}\Opt$ with $\Opt$ given by (\ref{final}).

{
\begin{figure}[htb!]
\begin{center}{\small
\begin{tabular}{c}
\includegraphics[width=0.9\textwidth]{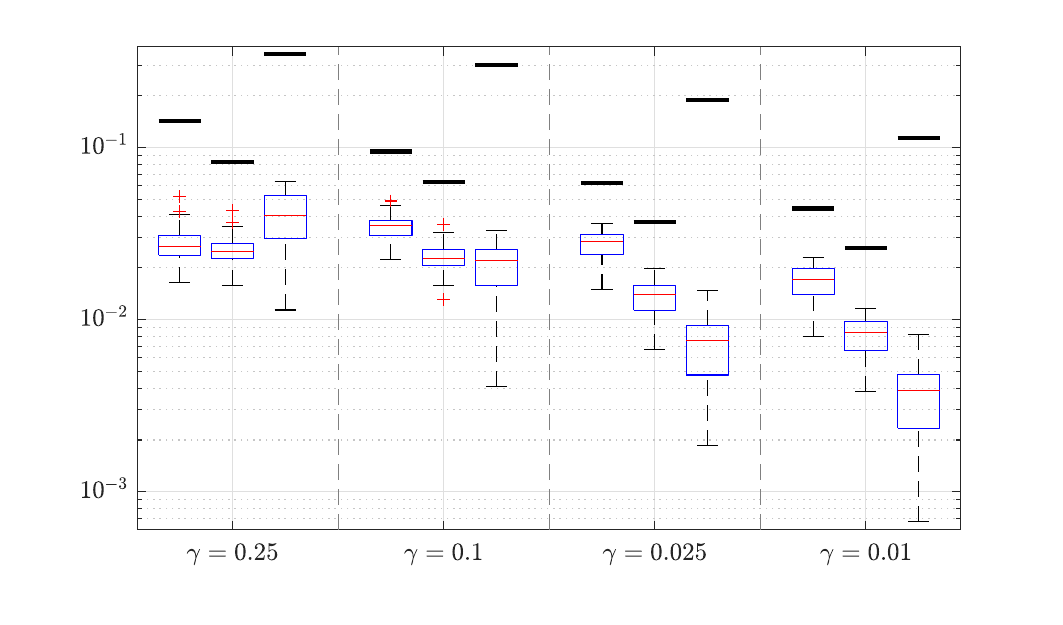}
\end{tabular}
\caption{\label{fig:2}  Distributions of $\ell_2$-recovery errors and upper bounds of the robust linear and robust polyhedral estimates for different values of $\gamma$ parameter.}
}\end{center}
\end{figure}
In Figure \ref{fig:2} we present error distributions and upper risk bounds (horizontal bar) of linear and polyhedral estimates in the numerical experiment with the model described in Section \ref{sectmod}.
In the plot cells, from left to right: (1) robust linear estimate by Proposition \ref{pr:lin_stoch} and upper bound $\mR$ on its $0.05$-risk, (2) robust linear estimate $w_1(\omega_1)$ yielded by Proposition \ref{proplin2} and upper bound $\wt\mR_1$ on its expected {error} risk, (3) robust polyhedral estimate by Proposition \ref{prop2} and upper bound on its $0.05$-risk.}
\subsubsection{A modification}
So far, our considerations related to polyhedral estimates were restricted to the case of sub-Gaussian $\eta$ and $\xi$. {Similarly to what was done in Section \ref{sectmod}, we} are about to show that passing from observation (\ref{eq1}) to its $K$-repeated, with  ``moderate'' $K$,  version (cf. \rf{eq1rep})
\[
\omega^K=\{\omega_k=A[\eta_k]x+ \xi_k,\;\;k=1,...,K\}
\]
with pairs $(\eta_k,\xi_k)$ independent across $k$, we can relax the sub-Gaussianity assumption replacing it with moment condition
\rf{relax}. Specifically, let us set
$$
\cH=\left\{h\in\bR^m: \sigma\|h\|_2\leq \tfrac{1}{8},\|\cA[h]\|_{\cX,2}\leq\tfrac{1}{8}\right\},\,\cA[h]x=[h^TA_1;...h^TA_q]
$$
(cf. (\ref{conda}) and (\ref{condb})).

Given tolerance an $m\times M$ contrast matrix $H$ with columns $h_j\in\cH$, and  observation (\ref{eq1rep}), we build the polyhedral estimate as follows.\footnote{Readers acquainted with the literature on robust estimation will immediately recognize that the proposed construction is nothing but a reformulation of the celebrated ``median-of-means'' estimate of \cite{nemirovskii1979complexity} (see also \cite{lerasle2011robust,hsu2014heavy,minsker2015geometric,lecue2020robust}) for our purposes. }
\begin{enumerate}
\item For $j=1,...,M$ we compute empirical medians $y_j$ of the data $h_j^T\omega_k$, $k=1,...,K$,
\[y_j=\med\{h_j^T\omega_k,\,1\leq k\leq K\}.
\]
\item We specify $\wh x^H(\omega^K)$ as a point from $\Argmin_{u\in \cX}\|y-H^TAu\|_\infty$ and use, as the estimate of $Bx$, the vector $\wh w^H_\poly(\omega^K)=B\wh x^H(\omega^K)$.
\end{enumerate}
\begin{lemma}\label{lem:mom} In the situation of this section, let $\xi_k$ and $\eta_k$ satisfy moment constraint of \rf{relax}, and let  $K\geq \ov \kappa= 2.5\ln[M/\epsilon]$.
Then estimate  $\wh w^H_\poly(\omega^K)$ satisfies
\[\Risk_{\epsilon}
[\wh w^H_\poly(\omega^K)|\cX] \leq \mP[H]
\]
 (cf. (\ref{easily})).
\end{lemma}
As an immediate consequence of the result of Lemma \ref{lem:mom}, the constructions and results of Sections \ref{sect323}--\ref{sectimpl} apply,
with  $\chi(\delta)=8$ and $\cH$ in the role of $\cH_\delta$, to our present situation in which the sub-Gaussianity of $\xi,\eta$ is relaxed
to the second moment condition \rf{relax} and instead of single observation $\omega$, we have access to a ``short''---with $K$ logarithmic
in $M/\epsilon$---sample of $K$ independent realizations of $\omega$.
\section{Uncertain-but-bounded perturbations}\label{uncbbnd}
In this section we assume that perturbation vector $\eta$ in (\ref{eq1}) is deterministic and runs through a given {\sl uncertainty set} $\cU$, so that (\ref{eq1}) becomes
\beq\label{eq1new}
\omega=A[\eta]x+\xi,\,\,A[\eta]=A+D[\eta],
\eeq
where $D[\eta]$ is (homogeneous) linear matrix-valued function of perturbation $\eta$ running through $\cU$.
As about observation noise $\xi$, we still assume that its distribution $P_x$ (which may depend on $x$) satisfies \rf{xi-sub}, i.e., is sub-Gaussian with zero mean and sub-Gaussian matrix parameter $\sigma^2I_m$ for every $x\in\cX$.

In our present situation it  is natural to redefine the notion of the $\epsilon$-risk of an estimate $\omega\mapsto\widehat{x}(\omega)$: here we consider uniform over $x\in \cX$ and $\eta\in \cU$ $\epsilon$-risk
\[
\Risk_{\epsilon}[\widehat{w}|\cX]=\sup_{x\in\cX,\eta\in \cU}\inf\Big\{\rho:\Prob_{\xi\sim P_x}\{\|\widehat{w}(A[\eta]x+\xi)-Bx\|>\rho\}\leq\epsilon\Big\}.
\]
Besides this, we, as before, assume that
$$
\|y\|=\max\limits_{\ell\leq L} \sqrt{y^TR_\ell y}\eqno{[R_\ell\succeq0,{\sum}_\ell R_\ell\succ0]}
$$

\subsection{Design of presumably good linear estimate}\label{sec:unbblin}
{Observe that the error of the linear estimate $\wh w^H(\omega)=H^T\omega$ satisfies
 \be
\|\wh w(A[\eta]x+\xi)-Bx\|\leq  \|H^T\xi\|+\max\limits_{x\in \cX,\eta\in\cU}\left\|H^T{D[\eta]} x\right\|+
\max\limits_{x\in \cX}\|[B-H^TA]x\|
\ee{riskdecnew}
}
Similarly to what was done in Section \ref{secplgle}, design of a presumably good linear estimate $\widehat{x}_H(\omega)$ consists in minimizing  over $H$ the sum of tight efficiently computable upper bounds on the terms in the right-hand side of \rf{riskdecnew}. Recall that bounds on the first and the last term were already established in Section \ref{secplgle} (cf. \rf{rell} and \rf{1stbpro} in the proof of Proposition \ref{pr:lin_stoch}). What is missing is a tight upper bound on
\[\ms(H)= \max\limits_{x\in \cX,\eta\in\cU}\left\|H^T{D[\eta]} x\right\|.
\]
In the rest of this section we focus on building efficiently computable upper bound on $\ms(H)$ which is convex in $H$; the synthesis of the contrast $H$ is then conducted by minimizing with respect to $H$ the resulting upper bound on estimation risk.
\par
We assume from now on that $\cU$ is a convex compact set in certain $\bR^q$. In this case $\ms(H)$ is what in \cite{JuKoNe} was called the {\sl robust norm}
$$
\begin{array}{c}
\|\cZ[H]\|_\cX=\max\limits_{Z\in\cZ[H]}\|Z\|_\cX,\,\,
\|Z\|_\cX=\max\limits_{x\in\cX}\|Zx\|\\
\end{array}
$$
of the {\sl uncertain $\nu\times n$ matrix}
$$
\cZ[H]=\{Z=H^TD[\eta]:\eta\in\cU\},
$$
i.e., the maximum, over {\sl instances} $Z\in\cZ[H]$, of operator norms of the linear mappings $x\mapsto Zx$ induced by the norm with the unit ball $\cX$ on the argument space and the norm $\|\cdot\|$ on the image space.
\par
It is well known that  aside of a very restricted family of special cases, robust norms do not allow for efficient computation. We are about to list known to us generic cases when these norms admit efficiently computable upper bounds which are tight within logarithmic factors.
\subsubsection{Scenario uncertainty}
This is the case where the nuisance set $\cU=\Conv\{\eta^1,...,\eta^S\}$ is given as a convex hull of moderate number of scenarios $\eta^s$.  In this case, $\ms(H)$ the maximum of operator norms:
\[
\ms(H)=\max\limits_{s\leq S}\max\limits_{x\in\cX}\|H^TD[\eta^s]x\|=\max\limits_{s\leq S,\ell\leq L} \|\cM_{s\ell}[H]\|_{\cX,2},\quad\cM_{s\ell}[H]=R_\ell^{1/2}H^TD[\eta^s],
\]
where, for $Q\in\bR^{\nu\times n}$,  $\|Q\|_{\cX,2}=\max\limits_{x\in\cX}\|Qx\|_2$ is the operator norm of the linear mapping $x\mapsto Qx:\bR^n\to\bR^\nu$ induced by the norm $\|\cdot\|_\cX$ with the unit ball $\cX$ on the argument space, and the Euclidean norm $\|\cdot\|_2$ on the image space. Note that this norm is efficiently computable in the  {\sl ellipsoid case} where $\cX=\{x\in\bR^n: x^TTx\leq1\}$ with $T\succ0$ (that is, for $K=1$, $T_1=T$, $\cT=[0,1]$ in (\ref{cX}))---one has
$
\|Q\|_{\cX,2}=\|Q{T^{-1/2}}\|_{2,2}.
$ When $\cX$ is a general ellitope,  norm $\|\cdot\|_{\cX,2}$ is difficult to compute. However,  it admits a tight efficiently computable convex in $Q$ upper bound:\footnote{We have already used it in the proof of Proposition \ref{pr:lin_stoch} when upper-bounding the corresponding terms $s_\ell(H)$ in the case of random uncertainty.}  it is shown in  \cite[Theorem 3.1]{JuKoNe} that function
$$
\Opt[Q]=\min\limits_{\lambda,\mu}\left\{\lambda+\phi_\cT(\mu):\mu\geq0,\hbox{\footnotesize$\left[\begin{array}{c|c}\lambda I_\nu&{1\over 2}Q\cr\hline {1\over 2}Q^T&{\sum}_k\mu_kT_k\cr\end{array}\right]$}\succeq0\right\}
$$
satisfies
$
\|Q\|_{\cX,2}\leq\Opt[Q]\leq 2.4\sqrt{\ln(4K)}\|Q\|_{\cX,2}
$.
As a result, under the circumstances,
\begin{align*}
\ov{\ms}(H)&=\max\limits_{{s\leq S,\ell\leq L}}\Opt_{s\ell}[H],\\
\Opt_{s\ell}[H]&=\min\limits_{{\lambda_\ell},{\mu^\ell}}\left\{{\lambda_\ell}+\phi_\cT({\mu^\ell}):{\mu^\ell}\geq0,
\hbox{\footnotesize$\left[\begin{array}{c|c}{\lambda_\ell} I_\nu&{1\over 2}R_\ell^{1/2}H^TD[\eta^s]\cr\hline {1\over 2}D^T[\eta^s]HR_\ell^{1/2}&{\sum}_k{\mu^\ell}_kT_k\cr\end{array}\right]$}\succeq0\right\},
\end{align*}
is a tight within the factor $2.4\sqrt{\ln(4K)}$ efficiently computable convex in $H$ upper bound on $\ms(H)$.

\subsubsection{Box and structured norm-bounded uncertainty}\label{intunc}
In the case of {\em structured norm-bounded uncertainty}  function $D[\eta]$ in the model \rf{eq1new} is of the form
\begin{align}\label{snbu}
D[\eta]&={\sum}_{\alpha=1}^qP_\alpha^T \eta_\alpha Q_\alpha\quad[P_\alpha\in\bR^{p_\alpha\times m},Q_\alpha\in \bR^{q_\alpha\times n}],\nn
\cU&=\{\eta=(\eta_1,...,\eta_q)\}=\cU_1\times...\times \cU_q,\\
\cU_\alpha&=\left\{\begin{array}{lll}\{\eta_\alpha=\delta I_{p_\alpha}:|\delta|\leq1\}\subset\bR^{p_\alpha\times p_\alpha},q_\alpha=p_\alpha&,\alpha\leq q_\s,&\hbox{\footnotesize["scalar perturbation blocks"]}\\
\{\eta_\alpha\in\bR^{p_\alpha\times q_\alpha}:\|\eta_\alpha\|_{2,2}\leq1\}&,q_\s<\alpha\leq {q}.&\hbox{\footnotesize["general perturbation blocks"]}
\end{array}\right.\nonumber
\end{align}
The special case of (\ref{snbu}) where $q_\s=q$, that is,
$$
\cU=\{\eta\in\bR^q:\|\eta\|_\infty\leq1\} \,\&\ A[\eta]=A+D[\eta]= A+{\sum}_{\alpha=1}^q\eta_\alpha A_\alpha
$$
is referred to as {\sl box uncertainty}.
In this section we operate with structured norm-bounded uncertainty (\ref{snbu}), assuming w.l.o.g. that all $P_\alpha$ are nonzero. The main result here (for underlying rationale and proof, see Section \ref{strnrm}) is as follows:
\begin{proposition}\label{prop912} Let $\cX\subset\bR^n$ be an ellitope:
$
\cX=P\cY$, where
$$
\cY=\{y\in\bR^n: \exists t\in\cT:y^TT_ky\leq t_k,k\leq K\}
$$
is a basic ellitope. Given the data of structured norm-bounded uncertainty (\ref{snbu}), consider the efficiently computable convex function
\begin{align*}
\ov\ms(H)&=\max\limits_{\ell\leq L}\Opt_\ell(H),\\
\Opt_\ell(H)&=\min\limits_{\mu,\upsilon,\lambda,U_s,V_s,U^t,V^t}
\Big\{\half[\mu+\phi_\cT(\upsilon)]:\,\mu\geq0,\upsilon\geq0,\lambda\geq0\\
&\left.\begin{array}{l}
\left[\begin{array}{c|c}U_s&-A_{s\ell}[H]P\cr\hline -P^TA_{s\ell}^T[H]&V_s\cr\end{array}\right]\succeq0,\,s\leq q_\s,\;
\left[\begin{array}{c|c}U^t&-L_{t\ell}^T[H]\cr\hline -L_{t\ell}[H]&\lambda_tI_{p_{q_\s+t}}\cr\end{array}\right]\succeq0,\,t\leq q-q_\s\\
V^t-\lambda_tP^TR_t^TR_tP\succeq0,\,t\leq q-q_\s\\
\mu I_\nu-{\sum}_sU_s-{\sum}_tU^t\succeq0, \;\;{\sum}_k\upsilon_k T_k-{\sum}_sV_s-{\sum}_tV^t\succeq0
\end{array}\right\}
\end{align*}
where
\begin{align*}
A_{s\ell}[H]&=R_\ell^{1/2}H^TP_s^TQ_s,\,1\leq s\leq q_\s\\
L_{t\ell}[H]&=P_{q_\s+t}HR_\ell^{1/2},\,R_t=Q_{q_\s+t},\,1\leq t\leq q-q_\s.
\end{align*}
Then
$$
\ms(H)\leq\ov\ms(H)\leq \varkappa(K)\max[\vartheta(2\kappa),\pi/2]\ms(H),
$$
where $\kappa=\max\limits_{\alpha\leq q_\s}\min[p_\alpha,q_\alpha]$ ($\kappa=0$ when $q_\s=0$),
$$
\varkappa(K)=\left\{\begin{array}{ll}1,&K=1,\\
{5\over 2}\sqrt{\ln(2K)},&K>1,
\end{array}\right.
$$ and $\vartheta(k)$ is a universal function of integer $k\geq0$ specified in (\ref{theta}) such that
$$
\vartheta(0)=0,\;\vartheta(1)=1,\;\vartheta(2)=\pi/2,\;\vartheta(3)=1.7348...,\;\vartheta(4)=2,\;\;
\vartheta(k)\leq\half {\pi\sqrt{k}},\;\;k\geq1.
$$
\end{proposition}
Note that the ``box uncertainty'' version of Proposition \ref{prop912} was derived in \cite{JuKoNe}.

\subsubsection{Robust estimation of linear forms}
Until now, we imposed no restrictions on the matrix $B$. We are about to demonstrate that when we aim at recovering the value of a given linear form $b^Tx$ of signal $x\in\cX$, i.e., when $B$ is a row vector:
\beq\label{linf}
Bx=b^Tx\qquad[b\in\bR^n],
\eeq
we can handle much wider family of uncertainty sets $\cU$ than those considered so far. Specifically, assume on the top of (\ref{linf}) that $\cU$ is a spectratope:
\begin{equation}\label{cUnew}
\begin{array}{c}
\cU=\{\eta=Qv,v\in\cV\},\,\cV=\{v\in\bR^M: \exists s\in\cS: S_\ell^2[v]\preceq s_\ell I_{d_\ell},\,\ell\leq L\},\\
S_\ell[v]={\sum}_{i=1}^Mv_iS^{i\ell},\,S^{i\ell}\in\bS^{d_\ell}\\
\end{array}
\end{equation}
(as is the case, e.g., with structured norm-bounded uncertainty)
and let $\cX$ be a spectratope as well:
\begin{equation}\label{cXnew}
\begin{array}{c}
\cX=\{x=Py,y\in\cY\},\,\cY=\{y\in\bR^N: \exists t\in\cT: T_k^2[y]\preceq t_k I_{f_k},\,k\leq K\},\\
T_k[y]={\sum}_{j=1}^Ny_jT^{jk},\,T^{jk}\in\bS^{f_k}.\\
\end{array}
\end{equation}
The contrast matrix $H$ underlying a candidate linear estimate becomes a vector $h\in\bR^m$, the associated  linear estimate being $\widehat{w}_h(\omega)=h^T\omega$. In our present situation $\nu=1$ we lose nothing when setting $\|\cdot\|=|\cdot|$. Representing $D[\eta]$ as ${\sum}_{\alpha=1}^q\eta_\alpha A_\alpha$, we get
\[
\mr_b(h)=\max\limits_{x\in\cX,\eta\in\cU} \left|h^T
{{\sum}}_\alpha \eta_\alpha A_\alpha x\right|=\max\limits_{\eta\in\cU,x\in\cX}\eta^TA[h]x,\quad A[h]=[h^TA_1;...;h^TA_q].
\]
In other words, $\mr_b(h)$ is the operator norm $\|A[h]\|_{\cX,\cU_*}$ of the linear mapping $x\mapsto A[h]x$ induced by the norm $\|\cdot\|_\cX$ with the unit ball $\cX$ on the argument space and the norm with the unit ball $\cU_*$---the polar of the spectratope $\cU$---on the image space. {Denote
\begin{align*}\lambda[\Lambda]&=[\Tr(\Lambda_1);...;\Tr(\Lambda_K)],\;\;\Lambda_k\in\bS^{f_k},\\
\lambda[\Upsilon]&=[\Tr(\Upsilon_1);...;\Tr(\Upsilon_L)],\;\;\Upsilon_\ell\in\bS^{d_\ell},
\end{align*}
and for $Y\in\bS^{d_\ell}$ and $X\in\bS^{f_k}$
\[R_\ell^{+,*}[Y]=\left[\Tr(YR^{i\ell}R^{j\ell})\right]_{i,j\leq M},\quad T_k^{+,*}[X]=\left[\Tr(XT^{ik}T^{jk})\right]_{i,j\leq N}.
\]}
Invoking \cite[Theorem 7]{JuKoNe}, we arrive at
\begin{proposition}\label{proplin3} In the case of (\ref{cUnew}) and (\ref{cXnew}), efficiently computable convex function
\begin{equation}\label{snormbound}
\overline{\mr}_b(h)=\min\limits_{\Lambda,\Upsilon}\left\{\half[\phi_{\cT}(\lambda[\Lambda])+\phi_{\cS}
(\lambda[\Upsilon]):
\begin{array}{l}\Lambda=\{\Lambda_k
\in
\bS^{f_k}_+,k\leq K\}, \Upsilon=\{\Upsilon_\ell\in\bS^{d_\ell}_+,\ell\leq L\}\\
\left[\begin{array}{c|c}{{\sum}}_\ell R_\ell^{+,*}[\Upsilon_\ell]&Q^T A[h]P\cr \hline P^T A^T[h] Q&{{\sum}}_kT_k^{+,*}[\Lambda_k]\cr\end{array}\right]\succeq0\\
\end{array}\right\}
\end{equation}
is a reasonably tight upper bound on $\mr_b(h)$:
\[
\mr_b(h)\leq \overline{\mr}_b(h)\leq \overline{\varsigma}\left({{\sum}}_{k=1}^Kf_k\right)\;\;
{\overline{\varsigma}\left({{\sum}}_{\ell=1}^Ld_\ell\right)}\mr_b(h)
\]
where $\overline{\varsigma}(J)=\sqrt{2\ln(5J)}$.
\end{proposition}
\subsection{Design of the robust polyhedral estimate}\label{sec:polyunbb}
On a close inspection, the strategy for designing a presumably good polyhedral estimate developed in Section \ref{pgle} for the case of random uncertainty works in the case of uncertain-but-bounded perturbations $A[\eta]=A+\underbrace{{{\sum}}_\alpha \eta_\alpha A_\alpha}_{D[\eta]}$,  $\eta\in\cU$, provided that the constraints (\ref{condab}) on the allowed columns $h$ of the contrast matrices are replaced with the constraint
\begin{subequations}
\label{condbbb}
\begin{align}\label{condbbb.a}
&\Prob_{\xi}\{|h^T\xi|>1/2\}\leq\delta/2,\\
\label{condbbb.b}
&\left|{\sum}_{\alpha=1}^q [h^TA_\alpha x]\eta_\alpha\right|\leq1/2\,\,\forall (x\in\cX,\eta\in\cU).
\end{align}
\end{subequations}
Assuming that $\cU$ and $\cX$ are the spectratopes (\ref{cUnew}), (\ref{cXnew}) and invoking Proposition \ref{proplin3}, an efficiently verifiable sufficient condition for $h$ to satisfy  the constraints (\ref{condbbb}) is
\begin{equation} \label{suffcond}
\|h\|_2\leq 2\sigma\sqrt{2\ln(2/\delta)}\;\;\mbox{and}\;\;
\overline{\mr}_b(h)\leq 1/2
\end{equation}
(see (\ref{conda}), (\ref{snormbound})). It follows that in order to build an efficiently computable upper bound for the $\epsilon$-risk of a polyhedral estimate associated with a given $m\times ML$ contrast matrix $H=[H_1,..,H_L]$, $H_\ell\in\bR^{m\times M}$, it suffices to check whether the columns of $H$ satisfy constraints (\ref{suffcond}) with $\delta=\epsilon/ML$. If the answer is positive, one can upper-bound the risk utilizing  the following spectratopic version of Proposition \ref{propanal}:

\begin{proposition}\label{propanalspect}
In the situation of this section, let $\epsilon\in (0,1)$, and let $H=[H_1,...,H_L]$ be $m\times ML$ matrix with $L$ blocks $H_\ell\in\bR^{m\times M}$  such that all columns of $H$ satisfy (\ref{suffcond}) with $\delta=\epsilon/ML$.  Consider optimization problem
\begin{align}\label{optprobspect}
\mP_+[H]&=2\min\limits_{\lambda_\ell,\Upsilon^\ell,\upsilon^\ell,\rho}\Big\{
\rho:\;\upsilon^\ell\geq0,\,\Upsilon^\ell=\{\Upsilon^\ell_k\in\bS^{f_k}_+,k\leq K\},\,\ell\leq L\\
&\qquad\left.\begin{array}{l}
\lambda_\ell+\phi_\cT(\lambda[\Upsilon^\ell])+{\sum}_{j=1}^M\upsilon^\ell_j \leq\rho,\,\ell\leq L\\
\left[\begin{array}{c|c}\lambda_\ell I_\nu&{1\over 2}R_\ell^{1/2}BP\cr\hline {1\over 2}{P^T}B^TR_\ell^{1/2}&P^TA^TH_\ell\Diag\{\upsilon^\ell\}H_\ell^TAP+
{\sum}_kT_k^{+,*}[\Upsilon^\ell_k]\cr\end{array}\right]\succeq0,\,\ell\leq L
\end{array}\right\}\nonumber
\end{align}
where \[
\lambda[\Upsilon^\ell]=[\Tr(\Upsilon^\ell_1);...;\Tr(\Upsilon^\ell_K)],\;\;\mbox{and}\;\;T_k^{+,*}(V)=
\left[\Tr(VT^{ik}T^{jk})\right]_{1\leq i,j\leq N}\hbox{\ for\ }V\in\bS^{f_k}.
\]
Then
\[
\Risk_{\epsilon}[\widehat{w}^H|\cX]\leq \mP_+[H].
\]
\end{proposition}
\paragraph{Remarks.}
As it  was already explained, when taken together, Propositions \ref{proplin3} and \ref{propanalspect} allow to compute efficiently
an upper bound on the $\epsilon$-risk of the polyhedral estimate associated with a given $m\times ML$ contrast matrix $H$: when the columns of $H$ satisfy (\ref{suffcond}) with $\delta=\epsilon/ML$, this bound is $\mP_+[H]$, otherwise it is, say, $+\infty$.
The outlined methodology can be applied to {\sl any} pair of spectratopes $\cX$, $\cY$.
However, to design a presumably good polyhedral estimate, we need to optimize the risk bound obtained in $H$,
and this seems to be difficult because the bound, same as its ``random perturbation'' counterpart, is nonconvex in $H$.
At present, we know only one generic situation where the synthesis problem admits ``presumably good'' solution---the case where both $\cX$ and $\cU$ are ellipsoids. Applying appropriate one-to-one linear transformations to perturbation $\eta$ and signal $x$, the latter situation can be reduced to that with
\beq\label{caseof}
\cX=\{x\in\bR^n:\|x\|_2\leq1\}, \;\;\cU=\{\eta\in\bR^q:\|\eta\|_2\leq1\},
\eeq
which we assume till the end of this section.
In  this case (\ref{suffcond}) reduces to
\beq\label{equivto}
\|h\|_2\leq [2\sigma\sqrt{2\ln(2/\delta)}]^{-1}\;\;\mbox{and}\;\; \|\cA[h]\|_{2,2}\leq 1/2
\eeq
where the matrix $\cA[h]$ is given by (\ref{barA}).
Note that (\ref{equivto}) is nothing but the constraint (\ref{eq100}) where the ellitope $\cX$ is set to be the unit Euclidean ball (that is, when $K=1$, $T_1=I_n$, and
 $\cT=[0;1]$ in (\ref{cX})) and the right hand side $\chi^{-1}(\delta)$ in the constraint is replaced with 1/2. As a result, (\ref{condbbb})
 can be processed in the same fashion as constraints (\ref{conda}) and (the single-ellipsoid case of) (\ref{eq100}) were processed in
 Sections  \ref{sect323} and \ref{sectstra} to yield a computationally efficient scheme for building a presumably good,
 in the case of (\ref{caseof}), polyhedral estimate. This scheme is the same as that described at the end of Section \ref{sec:end_stoch}
 with just one difference: the quantity $\chi(\delta)$ in the first semidefinite constraint of
 (\ref{final}) and (\ref{finall}) should now be replaced with constant $2$.
 Denoting by $\Opt$ the optimal value of the modified in the way just explained problem (\ref{final}), the $\epsilon$-risk of the polyhedral estimate yielded
 by an optimal solution to the problem is upper-bounded by $2\sqrt{\varkappa}\Opt$, with $\varkappa$ given by (\ref{varkappa}).
\appendix
\section{Proofs for Section \ref{secplgle}}
\subsection{Preliminaries: concentration of quadratic forms of sub-Gaussian vectors}
For the reader's convenience, we recall in this section some essentially known bounds for deviations of quadratic forms of sub-Gaussian random vectors (cf., e.g., \cite{hsu2012tail,rudelson2013hanson,spokoiny2013sharp}).
\paragraph{1$^o$.}
Let $\xi$  be a $d$-dimensional normal vector, $\xi\sim \cN(\mu,\Sigma)$. For all $h\in \bR^d$ and
$G\in \bS^d$ such that $G\prec \Sigma^{-1}$ we have the well known relationship:
\begin{align}
\ln\left(\bE_{\xi}\left\{e^{h^T\xi+{1\over 2} \xi^TG\xi}\right\}\right)&=-\half \ln \Det(I-\Sigma^{1/2}G\Sigma^{1/2})\nn
&+h^T\mu+\half \mu^TG\mu +\half [G\mu+h]^T\Sigma^{1/2}(I-\Sigma^{1/2}G\Sigma^{1/2})^{-1}\Sigma^{1/2}[G\mu-h].
\label{gausq}
\end{align}
Now, suppose that $\eta\sim \SG(0,\Sigma)$ where $\Sigma\in \bS^d_+$, let also $h\in \bR^d$ and $S\in \bR^{d\times d}$ such that $S\Sigma S^T\prec I$. Then for $\xi\sim \cN(h,S^TS)$ one has
\begin{align*}
\bE_{\eta}\left\{e^{h^T\eta+{1\over 2} \eta^TS^TS\eta}\right\}&=
\bE_{\eta}\left\{\bE_{\xi}\left\{e^{\eta^T\xi}\right\}\right\}=
\bE_{\xi}\left\{\bE_{\eta}\left\{e^{\eta^T\xi}\right\}\right\}\leq \bE_{\xi}\left\{e^{{1\over 2} \xi^T\Sigma \xi}\right\},
\end{align*}
so that
\begin{align*}
\ln\left(\bE_{\eta}\left\{e^{h^T\eta+{1\over 2} \eta^TS^TS\eta}\right\}\right)&\leq\ln\left(\bE_{\xi}\left\{e^{{1\over 2} \xi^T\Sigma \xi}\right\}\right)\\&=
-\half \ln \Det(I-S\Sigma S^T)+\half h^T\Sigma h+\half h^T\Sigma S^T(I-S\Sigma S^T)^{-1}S\Sigma h\\&=
-\half \ln \Det(I-S\Sigma S^T)+\half h^T\Sigma^{1/2}(I-S\Sigma S^T)^{-1}\Sigma^{1/2}h.
\end{align*}
In particular, when $\zeta\sim \SG(0,I)$, one has
\[
\ln\left(\bE_{\zeta}\left\{e^{h^T\zeta+{1\over 2} \zeta^TG\zeta}\right\}\right)\leq
-\half \ln \Det(I-G)
+\half h^T(I-G)^{-1}h=:\Phi(h,G).
\]
Observe that $\Phi(h,G)$ is convex  and continuous in $h\in \bR^d$ and $0\preceq G\prec I$ on its domain.
Using the inequality
(cf. \cite[Lemma 1]{laurent2000adaptive})
\be
\forall v\in[0,1[\quad -\ln(1-v)\leq {v}+{v^2\over 2(1-v)},
\ee{lm}
we get
\[
\Phi(h,G)\leq \half \Tr[ G]+\four \Tr[G(I-G)^{-1}G]+\half h^T(I-G)^{-1}h=:\wt\Phi(h,G).
\]
Finally, using
\[\Tr[G(I-G)^{-1}G]\leq (1-\lambda_{\max}(G))^{-1}\Tr[G^2],\quad h^T(I-G)^{-1}h\leq (1-\lambda_{\max}(G))^{-1}h^Th,
\]
we arrive at
\[
\wt\Phi(h,G)\leq \half \Tr[ G]+\four (1-\lambda_{\max}(G))^{-1}(\Tr[G^2]+2\|h\|_2^2)=:\ov\Phi(h,G).
\]

\paragraph{2$^o$.}
In the above setting, let $Q\in \bS^d_+$, $\alpha>2\lambda_{\max}(Q)$, $G=2Q/\alpha$, and let $h=0$. By the Cramer argument we conclude that
\be
\Prob\left\{\zeta^TQ\zeta\geq \alpha[\Phi(2Q/\alpha)+\ln \epsilon^{-1}]\right\}\leq \epsilon
\ee{dev11}
where $\Phi(\cdot)=\Phi(0,\cdot)$. In particular,
\be
\Prob\left\{\zeta^TQ\zeta\geq \min_{\alpha>2\lambda_{\max}(Q)}\alpha[\Phi(2Q/\alpha)+\ln \epsilon^{-1}]\right\}\leq \epsilon
\ee{dev1111}
Clearly, similar bounds hold with $\Phi$ replaced with $\wt\Phi$ and $\ov\Phi$. For instance,
\[
\Prob\left\{\zeta^TQ\zeta\geq \alpha[\ov\Phi(2Q/\alpha)+\ln \epsilon^{-1}]\right\}\leq \epsilon,
\]
so, when choosing $\alpha=2\lambda_{\max}(Q)+\sqrt{\Tr(Q^2)\over \ln\epsilon^{-1}}$ we arrive at  the ``standard bound''
\be
\Prob\left\{\zeta^TQ\zeta\geq \Tr(Q)+2\|Q\|_{\mathrm{Fro}}\sqrt{\ln\epsilon^{-1}}+2\lambda_{\max}(Q)\ln\epsilon^{-1}\right\}\leq \epsilon.
\ee{dev_simple0}
\begin{corollary}\label{cor227} Let $\epsilon\in (0,1)$, $W_1,...,W_L$ be matrices from $\bS^d_+$, and let $\upsilon\sim\SG(0,V)$ be a $d$-dimensional sub-Gaussian random vector. Then
$$
\Prob\left\{\max\limits_{\ell\leq L}\upsilon^TW_\ell\upsilon\geq \left[1+\sqrt{2\ln(L/\epsilon)}\right]^2\max\limits_{\ell\leq L}\Tr(W_\ell V)\right\}\leq \epsilon.
$$
\end{corollary}
{\bf Proof.} Let  $R^2=\max\limits_{\ell\leq L}\Tr(W_\ell V)$. W.l.o.g. we may assume that $\upsilon=V^{1/2}\zeta$ where $\zeta\sim \SG(0,I)$. Let us fix $\ell\leq L$. Applying \rf{dev_simple0} with $Q=V^{1/2}W_\ell V^{1/2}$ and $\epsilon$ replaced with  $\epsilon/L$, when taking into account that $\upsilon^TW_\ell\upsilon=\zeta^TQ\zeta$ with
\[\lambda_{\max}(Q)\leq \|Q\|_{\mathrm{Fro}}\leq \Tr(Q)\leq R^2,
\] we get
 \[\Prob\left\{\upsilon^TW_\ell\upsilon\geq \left[1+\sqrt{2\ln(L/\epsilon)}\right]^2R^2\right\}\leq {\epsilon\over L},
 \]
 and the claim of the corollary follows.\qed
\subsection{Proof of Proposition \ref{pr:lin_stoch}}
\label{sec:proof31}
Let $H$ be a candidate contrast matrix.
\paragraph{1$^o$.}
Observe that
 \begin{align}
\|\wh w^H(\omega)-Bx\|\leq \|H^T\xi\|+\left\|H^T{{\sum}}_{\alpha=1}^q\eta_\alpha A_\alpha x\right\|+\|[B-H^TA]x\|.
\label{3_terms}
\end{align}
Clearly,
\[\|[B-H^TA]x\|\leq \max_{\ell\leq L}\left\{\max_{x\in\cX}x^T[B-H^TA]^TR_\ell[B-H^TA]x\right\}^{1/2},
\]
so that by Theorem \ref{2020Prop4.6},
\be
\forall x\in \cX\quad \|[B-H^TA]x\|\leq \max\limits_{\ell\leq L}\mr_\ell(H)
\ee{r3}
where
\[\mr^2_\ell(H)=\min_\upsilon \left\{\phi_\cT(\upsilon):\,\upsilon\geq 0,\,\left[\begin{array}{c|c}I_\nu&R^{1/2}_\ell[B-H^TA]\cr\hline[B-H^TA]^TR^{1/2}_\ell&{\sum}_{k}\upsilon_kT_k\end{array}\right]\succeq 0\right\}.
\]
Taking into account that $\sqrt{u}=\min_{\lambda\geq 0}\{\tfrac{u}{4\lambda}+\lambda\}$ for $u> 0$, we get
\[\mr_\ell(H)=\min_{\upsilon,\lambda} \left\{\lambda+{\phi_\cT(\upsilon)\over 4\lambda}:\,\upsilon\geq 0,\lambda\geq 0,\,\left[\begin{array}{c|c}I_\nu&R^{1/2}_\ell[B-H^TA]\cr\hline[B-H^TA]^TR^{1/2}_\ell&{\sum}_{k}\upsilon_kT_k\end{array}\right]\succeq 0\right\}.
\]
Setting $\mu=\upsilon/(4\lambda)$, by the homogeneity of $\phi_\cT(\cdot)$ we obtain
\begin{align}\label{rell}
\mr_\ell(H)=\min_{\mu,\lambda} \left\{\lambda+\phi_\cT(\mu):\,\mu\geq 0,\,\left[\begin{array}{c|c}\lambda I_\nu&\half R^{1/2}_\ell[B-H^TA]\cr\hline\half [B-H^TA]^TR^{1/2}_\ell&{\sum}_{k}\mu_kT_k\end{array}\right]\succeq 0\right\}.
\end{align}
\paragraph{2$^o$.} Next, by Corollary \ref{cor227} of the appendix,
\be
\Prob \left\{\|H^T\xi\|\geq [1+\sqrt{2\ln(2L/\epsilon)}]\sigma\max_{\ell\leq L}\sqrt{\Tr(HR_\ell H^T)}\right\}\leq \epsilon/2.
\ee{1stbpro}
Similarly, because
\[\left\|H^T{{\sum}}_{\alpha=1}^q\eta_\alpha A_\alpha x\right\|=\max_{\ell\leq L}\left\|
R^{1/2}_\ell H^T[A_1x,...,A_qx]\eta\right\|_2,
\]we conclude that for any $x\in \cX$
\[\quad\Prob \left\{\left\|H^T{{\sum}}_{\alpha=1}^q\eta_\alpha A_\alpha x\right\|
\geq [1+\sqrt{2\ln(2L/\epsilon)}]\max_{\ell\leq L} s_\ell(H)\right\}\leq \epsilon/2\]
where $s_\ell(H)=\left\{\max_{x\in \cX} x^T\left[
{\sum}_\alpha A_\alpha^THR_\ell H^TA_\alpha\right]x\right\}^{1/2}$.
Again, by Theorem \ref{2020Prop4.6}, $s_\ell(H)$ may be tightly upper-bounded by the quantity $\ov{\ms}_\ell(H)$ such that
\[\ov{\ms}^2_\ell(H)=\min_\upsilon \left\{\phi_\cT(\upsilon):\,\upsilon\geq 0,\,\left[\begin{array}{c|c}I_{\nu q}
&[R^{1/2}_\ell H^TA_1;...;R^{1/2}_\ell H^TA_q]\cr\hline[A_1^THR^{1/2}_\ell,...,A_q^THR^{1/2}_\ell]&{\sum}_{k}\upsilon_kT_k\end{array}\right]\succeq 0\right\}.
\]
Now, repeating the steps which led to \rf{rell} above, we conclude that
\begin{align}\label{r2}
\ov{\ms}_\ell(H)&=\min_{\mu',\lambda'} \Big\{\lambda'+\phi_\cT(\mu'):\,\mu'\geq 0,\nn
&\qquad\quad\left[\begin{array}{c|c}\lambda' I_{\nu q}&\half[ R^{1/2}_\ell H^TA_1;...; R^{1/2}_\ell H^TA_q]\cr\hline\half [A_1^TH R^{1/2}_\ell,...,A_q^TH R^{1/2}_\ell]&{\sum}_{k}\mu'_kT_k\end{array}\right]\succeq 0\Bigg\}.
\end{align}
\paragraph{3$^o$.} When substituting the above bounds into \rf{3_terms}, we conclude that for every feasible solution
$\lambda_\ell,\mu^\ell, \kappa^\ell,\varkappa^\ell,\rho,\varrho$ to problem (\ref{prop31}) associated with $H$, the $\epsilon$-risk of the linear estimate $\widehat{w}^H_\lin(\cdot)$ may be upper-bounded by the quantity
$$
[1+\sqrt{2\ln(2L/\epsilon)}]\left[\sigma\max_{\ell\leq L}\|HR_\ell^{1/2}\|_\Fro+\rho\right]+\varrho.\eqno{\hbox{\qed}}
$$
\color{black}
\section{Proofs for Section \ref{pgle}}\label{proofprop}
\subsection{Proof of Proposition \ref{propanal}}
All we need to prove is that if $\lambda_\ell,\mu^\ell,\upsilon^\ell,\rho$ is a feasible solution to the optimization
problem (\ref{optprob}), then the inequality
\beq\label{then}
\Risk_{\epsilon}[\widehat{w}^H_\poly|\cX]\leq 2\rho
\eeq
holds. Indeed, let us fix $x\in\cX$. Since the columns of $H$ belong to $\cH_\delta$, the $P_x$-probability of the event
$$
\cZ^c=\{\zeta:\|H^T\zeta\|_\infty>1\}\eqno{[\zeta={\sum}_\alpha\eta_\alpha A_\alpha x+\xi]}
$$
is at most $ML\delta=\epsilon$. Let us fix observation $\omega=Ax+\zeta$ with $\zeta$ belonging to the complement $\cZ$ of $\cZ^c$. Then
$$
\|H^T[\omega-Ax]\|_\infty=\|H^T\zeta\|_\infty\leq 1,
$$
implying that the optimal value in the optimization problem $\min_{u\in\cX}\|H^T[Au-\omega\|_\infty$ is at most $1$. Consequently, setting $\wh{x}=\wh{x}^H(\omega)$, we have $\wh{x}\in\cX$ and $\|H^T[A\wh{x}-\omega]\|_\infty\leq1$, see (\ref{Au}). We conclude that setting $z=\half[x-\overline{x}]$, we have
$$
\|H_\ell^TAz\|_\infty\leq 1,\ell\leq L
$$
with $ z\in\cX$, implying that $z^TT_kz\leq t_k$, $k\leq K$, for some $t\in\cT$.
Now let $u\in\bR^\nu$  with $\|u\|_2\leq1$. Semidefinite constraints in (\ref{optprob}) imply that
\begin{align*}
u^TR_\ell^{1/2}Bz&\leq u^T\lambda_\ell I_\nu u+z^T\left[A^TH_\ell\Diag\{\upsilon^\ell\}H_\ell^TA+{\sum}_k\mu^\ell_kT_k\right]z\\
&\leq \lambda_\ell u^Tu+{\sum}_j\upsilon^\ell_j[H^TAz]_j^2+{\sum}_k\mu^\ell_kt_k\\
&\leq \lambda_\ell+{\sum}_j\upsilon^\ell_j +\phi_\cT(\mu^\ell)\leq\rho
\end{align*}
(recall that $\|u\|_2\leq 1$, $\lambda_\ell\geq0,\mu^\ell\geq0,\upsilon^\ell\geq0$, $t\in\cT$, and $\|H^T_\ell Az\|_\infty\leq1$).
 We conclude that $u^TR^{1/2}_\ell Bz\leq \rho$, $\ell\leq L$, whenever $\|u\|_2\leq1$, i.e., $\|R^{1/2}_\ell[Bz]\|_2\leq\rho^2$. The latter relation holds true for all $\ell\leq L$, implying that $\|Bz\|\leq \rho$, that is, $\|Bx-\widehat{x}(\omega)\|=2\|Bz\|\leq2\rho$ whenever $\zeta\in\cZ$. \qed
\subsection{Proof of Proposition \ref{prop2}}\label{sec:pr2pr}
\paragraph{0$^o$.} We need the following technical result.
\begin{theorem}\label{Khintchin}
  {\rm \cite[Theorem 4.6.1]{Tropp112}} Let $Q_i\in \bS^n$, $1\leq i\leq I$, and let $\xi_i$, $i=1,...,I$, be independent Rademacher ($\pm1$ with probabilities $1/2$) or $\cN(0,1)$ random variables. Then for all $t\geq 0$ one has \[ \Prob\left\{\left\|{{\sum}}_{i=1}^I\xi_i Q_i\right\|\geq t\right\}\leq 2n\exp\left\{-{t^2\over 2v_Q}\right\} \] where $\|\cdot\|$ is the spectral norm, and $ v_Q=\left\|{\sum}_{i=1}^IQ_i^2\right\|. $ \end{theorem}
\paragraph{1$^o$.} Proof of (i). Let $\lambda_j\geq0$, $g_j\in\cH$, $j\leq M$, and $\Sigma={\sum}_j\lambda_jg_jg_j^T$. Then for every $j$ there exists $r^j\in\cR$ such that $S_i^2[g_j]\preceq [r^j]_i I_{d_i}$, $i\leq I$. Assuming ${\sum}_j\lambda_j>0$ and setting $\kappa_j=[{\sum}_j\lambda_j]^{-1}\lambda_j$ and $r={\sum}_j\kappa_jr^j\in \cR$, we have
$$
\cS_i\left[{\sum}_j\lambda_jg_jg_j^T\right]={\sum}_j\lambda_jS_i^2[g_j]\preceq{\sum}_j\lambda_j[r^j]_iI_{d_i}=\left[{\sum}_j\lambda_j\right] r_iI_{d_i},
$$
implying that $(\Sigma,{\sum}_j\lambda_j)\in \bK$. The latter inclusion is true as well when $\lambda=0$.
\paragraph{2$^o$.} Proof of (ii). Let $(\Sigma,\rho)\in\bM$, and let us prove that $\Sigma={\sum}_{j=1}^{N} \lambda_jg_jg_j^T$ with $g_j\in\cH$, $\lambda_j\geq0$, and ${\sum}_j\lambda_j\leq\varkappa\rho$. There is nothing to prove when $\rho=0$, since in this case $\Sigma=0$ due to $(\Sigma,0)\in\bK$ combined with (\ref{when1b}). Now let $\rho>0$, so that for some $r\in \cR$ we have
\beq\label{eqeqqq}
\cS_i[\Sigma]\preceq \rho r_i I_{d_i},\,i\leq I,
\eeq
let $Z=\Sigma^{1/2}$, and let $O$ be the orthonormal $N\times N$ matrix of $N$-point Discrete Cosine Transform, so that all entries in $O$ are in magnitude $\leq\sqrt{2/N}$. For a Rademacher random vector $\varsigma=[\varsigma_1;...;\varsigma_M]$ (i.e., with entries $\varsigma_i$ which are independent Rademacher random variables), let
$$
Z^\varsigma=Z\Diag\{\varsigma\}O.
$$
In this case, one has
$
Z^\varsigma[Z^\varsigma]^T\equiv\Sigma,
$
that is,
\[
{\sum}_{p=1}^N\Col_p[Z^\varsigma]\Col_p^T[Z^\varsigma]\equiv \Sigma.
\]
Recall that
$$
\Col_j[Z^\varsigma]={\sum}_p
\varsigma_pO_{pj}\Col_p[Z],
$$
and thus
$$
S_i[\Col_j[Z^\varsigma]]={\sum}_p\varsigma_pO_{pj}S_i[\Col_p[Z]].
$$
Now observe that
\begin{align*}
{\sum}_p\big(O_{pj}S_i[\Col_p[Z]]\big)^2&={\sum}_pO_{pj}^2S_i^2[\Col_p[Z]]={\sum}_pO_{pj}^2\cS_i[\Col_p[Z]\Col_p^T[Z]]\\
\hbox{[see (\ref{when1a})]\;\;\;}&\preceq {2\over N}{\sum}_p\cS_i[\Col_p[Z]\Col_p^T[Z]]\\
&={2\over N}\cS_i[{\sum}_p\Col_p[Z]\Col_p^T[Z]]={2\over N}\cS_i[\Sigma]
\preceq {
2
\over N}{\rho}s_i I_{d_i}
\end{align*}
due to \rf{eqeqqq}.
By the Noncommutative Khintchine Inequality we have
\beq
\forall {\gamma}>0: \Prob\left\{S_i^2[\Col_j[Z^\varsigma]]\preceq \gamma{
2
\over N}\rho s_i I_{d_i}\right\}\geq 1-2d_i\exp\{-\gamma/2\}
\eeq
Setting
$$
\gamma=2\ln(4DN),\,D=
{\sum}_i d_i,\quad g_j^\varsigma=\sqrt{{N\over 2\gamma\rho}}\Col_j[Z^\varsigma],\quad\lambda_j=
{2\gamma
\rho\over N},\,1\leq j\leq N,
$$
we conclude that event
\[\Xi=\big\{\varsigma:S_i^2[g_j^\varsigma]\preceq s_i I_{d_i},i\leq I,j\leq N\big\}\subset \big\{g_j^\varsigma\in\cH,\,j\leq N\big\}
\]
satisfies $\Prob(\Xi)\geq \half$, while
\[{\sum}_j\lambda_jg_j^\varsigma[g_j^\varsigma]^T={\sum}_j\Col_j[Z^\varsigma]
\Col_j^T[Z^\varsigma]\equiv \Sigma\;\;\mbox{and}\;\;{\sum}_j\lambda_j=\gamma\kappa\rho=2\gamma\rho=\varkappa\rho.
\]
Thus, with probability $\geq 1/2$ (whenever $\varsigma\in\Xi$), vectors $g_j=g^\varsigma_j$ and $\lambda_j$ meet the requirements in (ii). \qed
\par
Note that the proof of the proposition suggests an efficient randomized algorithm for generating the required $g_j$ and $\lambda_j$: we generate realizations of $\varsigma$ of a Rademacher random vector, compute the corresponding vectors $g_j^\varsigma$, and terminate when all of them happen to belong to $\cH$. The corresponding probability not to terminate in course of the first $k$ rounds of randomization is then $\leq 2^{-k}$.
\subsection{Proof of Lemma \ref{lem:mom}}
The proof of the lemma is given by the standard argument underlying median-of-means construction (cf. \cite[Section 6.5.3.4]{nemirovskii1979complexity}). For the sake of completeness, we reproduce it here.
\paragraph{1$^o$.} Observe that when \rf{relax} holds, $h\in\cH$, $x\in\cX$ and $\zeta=\xi+\sum_\alpha\eta_\alpha A_\alpha x$, the probability of the event
$$
\{|h^T\zeta|>1\}
$$
is at most 1/8.
Indeed, when $|h^T\zeta|>1$ implies that either $|h^T\xi|>1/2$ or $|\eta^T\cA[h]x|>1/2$. By the Chebyshev inequality, the probability of the first of these events is at most $4\bE\{(h^T\xi)^2\}\leq 4\sigma^2\|h\|_2^2\leq{1\over 16}$ (we have used the first relation in (\ref{relax}) and took into account that $h\in\cH$). By similar argument, the probability of the second event is at most $4\bE\{(\eta^T\cA[h]x)^2\}\leq 4\|\cA[h]x\|_2^2\leq {1\over 16}$.
\paragraph{2$^o$.}
Let $\zeta_k=\omega_k-Ax$. By construction, $z_j=y_j-h_j^TAx$ is the median of the i.i.d. sequence $h_j^T\zeta_k$, $k=1,...,K$. When $|z_j|>1$, at least $K/2$ of the events $\{|h_j^T\zeta_k|>1\}$, $k\leq K$, take place. Because the probability of each of $K$ independent events is $\leq 1/8$, it is easily seen\footnote{We refer to, e.g., \cite[Section 2.3.2]{EuclSep} for the precise justification of this obvious claim.} that the probability that at least $K/2$ of them happen is bounded with
\[
\pi(K):=\sum_{k\geq K/2}\binom{K}{k}(1/8)^k(7/8)^{K-k}\leq \sum_{k\geq K/2}\binom{K}{k}2^{-K}[(1/4)^k(7/4)^{K-k}]\leq(\sqrt{7}/4)^K\leq e^{-0.4K}.
\]
In other words, the probability of each event $E_j=\{\omega^K: |y_j-h_j^TAx|>1\}$, $j=1,...,M$, is bounded with $\pi(K)$. Thus, none of the events $E_1,...,E_M$ takes place with probability at least $1-M\pi(K)$, and in such case we have $\|y-H^TAx\|_\infty\leq1$, and so $\|y-H^TA\wh x^H(\omega^K)\|_\infty\leq1$ as well. We conclude that for every $x\in\cX$, the probability of the event
$$
\left\{x-\wh x^H(\omega^K)\in2\cX,\, \|H^TA[x-\wh x^H(\omega^K)]\|_\infty\leq 2\right\}
$$
is at least  $1-M\pi(K)\geq 1-\epsilon$ when $K\geq 2.5\ln[M/\epsilon]$, and when it happens, one has $\|Bx-\wh w^H_\poly(\omega^K)\|\leq\mP[H]$.\qed
\subsection{Proof of Proposition \ref{proplin2}}
\paragraph{1$^o$.} Let $\ell\leq L$ and $k\leq K$ be fixed, let $H=H_\ell\in\bR^{m\times \nu}$ be a candidate contrast matrix, and let $\lambda,\mu,\kappa,\varkappa$ be a feasible solution to \rf{add0}.  One has
\begin{align}
\bE_{\xi_k} \left\{\|R_\ell^{1/2}H^T\xi_k \|_2^2\right\}&=\Tr\left(\bE_{\xi_k} \left\{R_\ell^{1/2}H^T\xi_k \xi_k ^THR_\ell^{1/2}\right\}\right)
\label{add1}
\leq \sigma^2\Tr(HR_\ell H^T)=\sigma^2\|HR^{1/2}_\ell\|_\Fro^2.
\end{align}
Next, for  any $x\in \cX$  fixed we have
\begin{align}\label{add2}
\bE_{\eta_k} \left\{\left\|R_\ell^{1/2}H^T[{\sum}_\alpha[\eta_k]_\alpha A_\alpha]x\right\|_2^2\right\}&=
\bE_{\eta_k} \left\{\left\|
R^{1/2}_\ell H^T[A_1x,...,A_qx]\eta_k \right\|_2^2\right\}=x^T\left[
{\sum}_\alpha A_\alpha^THR_\ell H^TA_\alpha\right]x\nn
&=\|[R_\ell^{1/2}H^TA_1;...;R_\ell^{1/2}H^TA_q]x\|_2^2\leq(\lambda+\phi_\cT(\mu))^2
\end{align}
where the concluding inequality follows from the constraints in (\ref{add0}) (cf. item 2$^o$ of the proof of Proposition \ref{pr:lin_stoch}).
Next, similarly to item 1$^o$ of the proof of Proposition \ref{pr:lin_stoch} we have
$$
\|R_\ell^{1/2}(B-H^TA)x\|_2^2\leq (\kappa+\phi_\cT(\varkappa))^2.
$$
Put together, the latter bound along with \rf{add1} and \rf{add2} imply \rf{expkl}.
\paragraph{2$^o$.} By the Chebyshev inequality,
\[
\forall \ell,k\quad \Prob\left\{\|R_\ell^{1/2}(w_{\ell}(\omega_k)-Bx)\|_2\geq 2\wt\mR_\ell[H_\ell]\right\}\leq \four;
\]
applying  \cite[Theorem 3.1]{minsker2015geometric} 
we conclude that
\[
\forall \ell\quad \Prob\left\{\|R_\ell^{1/2}(z_{\ell}(\omega^K)-Bx)\|_2\geq 2C_\alpha \wt\mR_\ell[H_\ell]\right\}\leq e^{-K\psi(\alpha,{1\over 4})}
\]
where
\be
\psi(\alpha,\beta)=(1-\alpha)\ln {1-\alpha\over 1-\beta}+\alpha \ln{\alpha\over \beta}
\ee{psidef}
and $C_\alpha={1-\alpha\over \sqrt{1-2\alpha}}$. When choosing $\alpha={\sqrt{3}\over 2+\sqrt{3}}$ which corresponds to $C_\alpha=2$ we obtain $\psi(\alpha,\four)=0.1070...$ so that for $\ell\leq L$
\[
 \Prob\left\{\|R_\ell^{1/2}(z_{\ell}(\omega^K)-Bx)\|_2\geq  4 {\wt{\mR}}_\ell[H_\ell]\right\}\leq e^{-0.1070 K}
 \]
 what is \rf{prbl}.
 \paragraph{3$^o$.} Now, let $K\geq\ln(L/\epsilon)/0.1070$. In this case, for all $\ell\leq L$
\[
 \Prob\left\{\|R_\ell^{1/2}(z_{\ell}(\omega^K)-Bx)\|_2\geq  4{\wt{\mR}}_\ell[H_\ell]\right\}\leq \epsilon/L,
\]
so that with probability $\geq 1-\epsilon$ the set $ \cW(\omega^K)$ is not empty (it contains $Bx$), and for all $v\in \cW(\omega^K)$ one has
\[\|R_\ell^{1/2}(v-Bx)\|_2\leq \|R_\ell^{1/2}(z_{\ell}(\omega^K)-v)\|_2+\|R_\ell^{1/2}(z_{\ell}(\omega^K)-Bx)\|_2\geq  8 {\wt{\mR}}_\ell[H_\ell].\eqno{\mbox{\qed}}
\]
{
\section{Proofs for Section \ref{uncbbnd}}}
\subsection{Proof of Proposition  \ref{propanalspect}}
The proof follows that of Proposition \ref{propanal}.
All we need to prove is that if $H$ satisfies the premise of the proposition and $\lambda_\ell,\Upsilon^\ell,\upsilon^\ell,\rho$ is a feasible solution to (\ref{optprobspect}), then the inequality
\beq\label{thenspect}
\Risk_{\epsilon}[\widehat{w}^H_\poly|\cX]\leq 2\rho
\eeq
holds. Indeed, let us fix $x\in\cX$ and $\eta\in\cU$. Since the columns of $H$ satisfy (\ref{suffcond}),
 the $P_x$-probability of the event
$$
\cZ_{x,\eta}=\{\xi:\|H^T[D[\eta]x+\xi\|_\infty\leq 1\}
$$
is at least $1-ML\delta=1-\epsilon$. Let us fix observation $\omega=Ax+D[\eta]x+\xi$ with $\xi\in\cZ_{x,\eta}$. Then
\begin{equation}\label{added}
\|H^T[\omega-Ax]\|_\infty=\|H^T[D[\eta]x+\xi]\|_\infty\leq 1,
\end{equation}
implying that the optimal value in the optimization problem $\min_{u\in\cX}\|H^T[Au-\omega\|_\infty$ is at most $1$.  Consequently, setting $\wh{x}=\wh{x}^H(\omega)$, we have $\wh{x}\in\cX$ and $\|H^T[A\wh{x}-\omega]\|_\infty\leq1$, see (\ref{Au}). These observations combine with (\ref{added}) and the inclusion $x\in\cX$ to imply that for $z=\half[x-\wh x]$ we have $z\in\cX$ and
$\|H^Tz\|_\infty\leq1$. Recalling what $\cX$ is we conclude that $z=Py$ with $T_k^2[y]\preceq t_kI_{f_k},k\leq K$ for some $t\in\cT$ and
\begin{equation}\label{werwe}
\|H_\ell^TAPy\|_\infty=\|H_\ell^TAz\|_\infty\leq 1,\;\ell\leq L.
\end{equation}
Now let $u\in\bR^\nu$  with $\|u\|_2\leq1$. Semidefinite constraints in (\ref{optprobspect}) imply that
\begin{align}\label{comp}
u^TR_\ell^{1/2}Bz&=u^TR_\ell^{1/2}BPy\leq u^T\lambda_\ell I_\nu u+y^T\left[PA^TH_\ell\Diag\{\upsilon^\ell\}H_\ell^TAP+{\sum}_kT_k^{+,*}[\Upsilon^\ell_k]\right]y\nn
&= \lambda_\ell u^Tu+{\sum}_j\upsilon^\ell_j\underbrace{[H_\ell^TAPy]_j^2}_{\leq1\hbox{ by
\rf{werwe}}}+{\sum}_ky^TT_k^{+,*}[\Upsilon^\ell_k]y\nn
&\leq\lambda_\ell+{\sum}_j\upsilon^\ell_j +\sum_k\sum_{i,j\leq N}y_iy_j\Tr(\Upsilon^\ell_kT^{ik}T^{jk})\nn
&=\lambda_\ell+{\sum}_j\upsilon^\ell_j +\sum_k\Tr(\Upsilon^\ell_kT^2_k[y])\nn
&\leq \lambda_\ell+{\sum}_j\upsilon^\ell_j +\sum_kt_k\Tr(\Upsilon^\ell_k) \hbox{\ [due to $\Upsilon^\ell\succeq0$ and $T_k^2[y]\preceq t_kI_{f_k}$]}\nn
&\leq \lambda_\ell+{\sum}_j\upsilon^\ell_j+\phi_\cT(\lambda[\Upsilon^\ell])\leq\rho
\end{align}
where the concluding inequality follows from the constraints of \rf{optprobspect}.
(\ref{comp}) holds true for all $u$ with ${\|u\|_2}\leq1$, and we conclude that  for $x\in\cX$ and $\eta\in\cU$  and $\xi\in\cZ_{x,\eta}$ (recall that the latter inclusion takes place with $P_x$-probability $\geq1-\epsilon$) we have
\[\|R_\ell^{1/2}B[\wh x^H(Ax+D[\eta]
x+\xi)-x]\|_2\leq 2\rho, \;\ell\leq L.
\] Recalling what $\|\cdot\|$ is, we get
$$
\forall (x\in\cX,\eta\in\cU): \Prob_{\xi\sim P_x}\{\|B[x-\wh x^H(Ax+D[\eta]x+\xi)]\|>{2}\rho\|\leq\epsilon,
$$
that is, $\Risk_\epsilon[\wh w^H_\poly|\cX]\leq{2}\rho$. The latter relation holds true whenever $\rho$ can be extended to a feasible solution to (\ref{optprobspect}), and (\ref{thenspect}) follows. \qed

\subsection{Robust norm of uncertain matrix with structured norm-bounded uncertainty}\label{strnrm}
\subsubsection{Situation and goal}
Let matrices $A_s\in\bR^{m\times n}$, $s\leq S$, and $L_t\in\bR^{p_t\times m}$, $R_t\in\bR^{q_t\times n}$, $t\leq T$, be given. These data specify {\sl uncertain $m\times n$ matrix}
\begin{equation}\label{cA}
\cA=\{A={\sum}_s\delta_sA_s+{\sum}_tL_t^T\Delta_tR_t:|\delta_s|\leq1\,\forall s\leq S, \|\Delta_t\|_{2,2}\leq1\,\forall t\leq T\}.
\end{equation}
Given ellitopes
\begin{equation}\label{ellitopes}
\begin{array}{rcl}
\cX&=&\{Py: y\in\cY\}\subset\bR^n, \,\cY=\{y\in\bR^N\ \&\ \exists t\in\cT: y^TT_ky\leq t_k,k\leq K\},\\
\cB_*&=&\{Qz:z\in\cZ\}\subset\bR^m,\,\cZ=\{z\in\bR^M:\exists s\in\cS: z^TS_\ell z\leq s_\ell,\,\ell\leq L\},
\end{array}
\end{equation}
we want to upper-bound the robust norm
$$
\|\cA\|_{\cX,\cB}=\max_{A\in\cA}\|A\|_{\cX,\cB},
$$
of uncertain matrix $\cA$ induced by the norm $\|\cdot\|_\cX$ with the unit ball $\cX$ in the argument  space and the norm $\|\cdot\|_\cB$ with the unit ball $\cB$ which is the polar of $\cB_*$ in the image space.
\subsubsection{Main result}
\begin{proposition}\label{propmain} Given uncertain matrix (\ref{cA}) and ellitopes (\ref{ellitopes}), consider convex optimization problem
\begin{subequations}\label{cprb}
\begin{align}
\Opt&=\min\limits_{\mu,\upsilon,\lambda,\atop U_s,V_s,U^t,V^t}\half [\phi_\cS(\mu)+\phi_\cT(\upsilon)]\nn
&\mathrm{subject\;to}\nn
\label{cprb.a}
&\mu\geq0,\,\upsilon\geq0,\,\lambda\geq0\nn
&\left[\begin{array}{c|c}U_s&-Q^TA_sP\cr\hline -P^TA_s^TQ&V_s\end{array}\right]\succeq0\\
\label{cprb.b}&\left[\begin{array}{c|c}U^t&-Q^TL_t^T\cr\hline -L_tQ&\lambda_tI_{p_t}\end{array}\right]\succeq0,\;
V^t-\lambda_tP^TR_t^TR_tP\succeq0\\
\label{cprb.c}&{\sum}_\ell\mu_\ell S_\ell-{\sum}_sU_s-{\sum}_tU^t\succeq0\\
\label{cprb.d}&{\sum}_k\upsilon_k T_k-{\sum}_sV_s-{\sum}_tV^t\succeq0
\end{align}
\end{subequations}

The problem is strictly feasible and solvable, and
\begin{equation}\label{and}
\|\cA\|_{\cX,\cB}\leq\Opt\leq \varkappa(K)\varkappa(L)
\max\left[\vartheta(2\kappa),\pi/2\right]\|\cA\|_{\cX,\cB}
\end{equation}
where
\begin{itemize}
\item the function $\vartheta(k)$ of nonnegative integer $k$ is given by $\vartheta(0)=0$ and
{\small\begin{equation}\label{theta}\vartheta(k)=\left[\min_{\alpha}\left\{(2\pi)^{-k/2}\int
|\alpha_1u_1^2+...+\alpha_ku_k^2|\e^{-{u^T u/2}}du,\;\alpha\in{\bR}^k,\|\alpha\|_1=1\right\}\right]^{-1},\;\;k\geq1;
\end{equation}}
\item $\kappa=\max\limits_{s\leq S}\rank(A_s)$ when $S\geq1$, otherwise $\kappa=0$;
\item $\varkappa(\cdot)$ is given by
\begin{equation}\label{varkappaJ}
\varkappa(J)=\left\{\begin{array}{ll}1,&J=1,\\
{5\over 2}\sqrt{\ln(2J)},&J>1.
\end{array}\right.
\end{equation}
\end{itemize}
\end{proposition}
\paragraph{Remarks.}
The rationale behind (\ref{cprb}) is as follows. Checking that the $\cX,\cB$-norm of uncertain $m\times n$ matrix (\ref{cA}) is $\leq a\in\bR$  is the same as to verify that for all $\delta_s\in[-1,1],\;\Delta_t:\|\Delta_t\|_{2,2}\leq1$
\[
{\sum}_s\delta_su^TA_sv+{\sum}_tu^TL_t^T\Delta_tR_tv\leq a\|u\|_{\cB_*}\|v\|_\cX\quad\forall (u\in\bR^m,v\in\bR^n),
\]
or, which is the same due to what $\cB_*$ and $\cX$ are, that for all $\delta_s\in[-1,1],\Delta_t:\|\Delta_t\|_{2,2}\leq1$
\be
{\sum}_s\delta_sz^TQ^TA_sPy+{\sum}_tz^TQ^TL_t^T\Delta_tR_tPy\leq a\|z\|_{\cZ}\|y\|_\cY\quad\forall (z\in\bR^M,y\in\bR^N).
\ee{rela}
A simple certificate for (\ref{rela}) is a collection of positive semidefinite matrices $U_s,V_s,U^t,V^t$, $U,V$  such that for all $z\in\bR^M,\,y\in\bR^N$ and all $s\leq S$, $t\leq T$ it holds
\begin{subequations}
\label{suchsuch}
\begin{align}
2z^T[Q^TA_sP]y&\leq z^TU_sz+y^TV_sy,\label{suchsuch.a}\\
2z^TQ^TL_t^T\Delta_tR_tPy&\leq z^TU^tz+y^TV^ty\;\;\forall (\Delta_t:\|\Delta_t\|_{2,2}\leq1),\label{suchsuch.b}\\
{\sum}_sU_s+{\sum}_tU^t&\preceq U,\label{suchsuch.c}\\
{\sum}_sV_s+{\sum}_tV^t&\preceq V,\label{suchsuch.d}\\
\max\limits_{z\in\cZ} z^TUz+\max\limits_{y\in\cY}y^TVy&\leq 2a.\label{suchsuch.e}
\end{align}\end{subequations}
Now, \rf{suchsuch.a} clearly is the same as \rf{cprb.a}. It is known (this fact originates from \cite{LMI}) that \rf{suchsuch.b} is the same as existence of $\lambda_t\geq0$ such that \rf{cprb.b} holds. Finally, existence of $\mu\geq0$ such that
${\sum}_\ell\mu_\ell S_\ell\succeq U$ and $\upsilon\geq0$ such that
${\sum}_k\upsilon_k T_k\succeq V$ (see \rf{cprb.c} and \rf{cprb.d}) implies due to the structure of $\cZ$ and $\cY$ that $\max_{z\in\cZ}z^TUz\leq\phi_\cS(\mu)$ and $\max_{y\in\cY}y^TVy\leq\phi_\cT(\upsilon)$. The bottom line is that a feasible solution to (\ref{cprb}) implies the existence of a certificate \[\left\{U_s,U^t,V_s,V^t,s\leq S,t\leq T,U={\sum}_\ell\mu_\ell S_\ell,V={\sum}_k\upsilon_kT_k\right\}\] for relation
(\ref{rela}) with $a=\half[\phi_\cS(\mu)+\phi_\cT(\upsilon)]$.
\par\noindent
{\bf Proof of Proposition \ref{propmain}.} {\bf 1}$^o$. Strict feasibility and solvability of the problem are immediate consequences of ${\sum}_\ell S_\ell\succ0$ and ${\sum}_kT_k\succ0$.
\par Let us prove the first inequality in (\ref{and}). All we need to show is that if \begin{itemize}
\item[][a] $\mu,\upsilon,\lambda,U_s,V_s,U^t,V^t$ is feasible for (\ref{cprb}),
\item[][b]
$x=Py$ with $y^TT_ky\leq \tau_k$, $k\leq K$,  for some $\tau\in\cT$ and $u=Qz$ for some $z$ such that $z^TS_\ell z\leq \varsigma_\ell$, $\ell\leq   L$, for some $\varsigma\in\cS$, and
\item[][c] $\delta_s$, $\Delta_t$ satisfy $|\delta_s|\leq1$, $\|\Delta_t\|_{2,2}\leq1$,
\end{itemize}
then  $\gamma:=u^T[{\sum}_s\delta_sA_s+{\sum}_tL_t^T\Delta_t R_t]x\leq {1\over 2}[\phi_\cS(\mu)+\phi_\cT(\upsilon)]$. Assuming [a--c], we have
\begin{align*}
\gamma&={\sum}_s\delta_s z^TQ^TA_sPy+{\sum}_tz^TQ^TL_t^T\underbrace{\Delta_tR_tPy}_{\zeta_t}\\
&\leq \half z^T\left[{\sum}_sU_s\right]z+\half y^T\left[{\sum}_s V_s\right]y+{\sum}_t\|L_tQz\|_2\|\zeta_t\|_2\quad
\hbox{\ [by \rf{cprb.a} and due to $|\delta_s|\leq1$]}\\
&\leq \half z^T\left[\sum_s U_s\right]z+\half y^T\left[\sum_s V_s\right]y+\sum_t\sqrt{(\lambda_tz^TU^tz)(y^TP^TR_t^TR_tPy)}\\&\quad\hbox{[due to \rf{cprb.b} and $\|\Delta_t\|_{2,2}\leq1$]}\\
&=\half z^T\left[{\sum}_sU_s\right]z+\half y^T\left[{\sum}_s V_s\right]y+{\sum}_t\sqrt{(z^TU^tz)(\lambda_ty^TP^TR_t^TR_tPy)}. \end{align*}
Thus, by the second inequality of \rf{cprb.b},
\begin{align}\label{chain74}
\gamma&\leq \half z^T\left[{\sum}_sU_s\right]z+\half y^T\left[{\sum}_s V_s\right]y+{\sum}_t\sqrt{(z^TU^tz)(y^TV^ty)}\nn
&\leq \half z^T\left[{\sum}_sU_s\right]z+\half y^T\left[{\sum}_s V_s\right]y+\half{\sum}_t[z^TU^tz+y^TV^ty]\nn
&= \half\left[z^T\left[{\sum}_sU_s+{\sum}_tU^t\right]z+y^T\left[{\sum}_sV_s+{\sum}_tV^t\right]y\right]\nn
&\leq \half\left[{\sum}_\ell\mu_\ell z^TS_\ell z+{\sum}_ky^T\upsilon_k T_ky\right]\quad \hbox{[by \rf{cprb.c} and \rf{cprb.d}]}\nn
&\leq  \half\left[{\sum}_\ell\mu_\ell \varsigma_\ell+{\sum}_k\upsilon_k\tau_k\right]\quad \hbox{[due to $z^TS_\ell z\leq \varsigma_\ell,
y^TT_ky\leq \tau_k$]}\nn
&\leq \half[\phi_\cS(\mu)+\phi_\cT(\upsilon)]\quad \hbox{[since $\varsigma\in\cS,\tau\in\cT$]}
\end{align}
as claimed.
\paragraph{2$^o$.} Now, let us prove the second inequality in (\ref{and}). Observe that
$$
 \bS=\{0\}\cup\{[s;\sigma]:\sigma>0,s/\sigma\in\cS\},\quad \bT=\{0\}\cup\{[t;\tau]:\tau>0,t/\tau\in\cT\},
$$
are regular cones with the duals
$$
\bS_*=\{[g;\sigma]:\sigma\geq\phi_\cS(-g)\},\quad\bT_*=\{[h;\tau]:\tau\geq\phi_\cT(-h)\},
$$
and (\ref{cprb}) can be rewritten as the conic problem
\begin{align}
2\Opt&=\min\limits_{{\alpha,\beta,\mu,
\upsilon,\atop\lambda,U_s,V_s,U^t,V^t}}\alpha+\beta\tag{P}\\
&\hbox{subject to}\nn
&[-\mu;\alpha]^{\hbox{\footnotesize$[\overline{g},\overline{\alpha}]$}}\in\bS_*,\,[-\upsilon;\beta]^{\hbox{\footnotesize$[\overline{h},\overline{\beta}]$}}\in\bT_*,
\,\mu^{\hbox{\footnotesize$\overline{\mu}$}}\geq0,\,\upsilon^{\hbox{\footnotesize$\overline{\upsilon}$}}\geq0,\,\lambda^{\hbox{\footnotesize$\overline{\lambda}$}}\geq0\nn
&\left[\begin{array}{c|c}U_s&-Q^TA_sP\cr\hline -P^TA_s^TQ&V_s\cr\end{array}\right]^{\hbox{\footnotesize$\left[\begin{array}{c|c}\overline{U}_s&\overline{A}_s\cr\hline
\overline{A}_s^T&\overline{V}_s\end{array}\right]$}}
\succeq0, \;s\leq S\nn
&\left[\begin{array}{c|c}U^t&-Q^TL_t^T\cr\hline -L_tQ&\lambda_tI_{p_t}\cr\end{array}\right]^{\hbox{\footnotesize$\left[\begin{array}{c|c}\overline{U}^t&\overline{L}_t^T\cr\hline
\overline{L}_t&\overline{\Lambda}_t\cr\end{array}\right]$}}\succeq0,\;\;{[V^t-\lambda_tP^TR_t^TR_tP]}^{\hbox{\footnotesize$\overline{V}^t$}}\succeq0,\;t\leq T\nn
&{[{\sum}_\ell\mu_\ell S_\ell-{\sum}_sU_s-{\sum}_tU^t]}^{\hbox{\footnotesize$\overline{S}$}}\succeq0, \;\;{[{\sum}_k\upsilon_k T_k-{\sum}_sV_s-{\sum}_tV^t]}^{\hbox{\footnotesize$\overline{T}$}}\succeq0\label{tagP.ef}
\end{align}
(superscripts are the Lagrange multipliers for the corresponding constraints).
$(P)$ clearly is solvable and strictly feasible, so that $2\Opt$ is the optimal value of the (solvable!) conic dual of $(P)$:

\begin{align}\tag{D}
2\Opt&=\max\limits_{{\overline{\alpha},\overline{\beta},\overline{g},\overline{h},\overline{\mu},\overline{\upsilon},\overline{\lambda},\overline{S},\overline{T},\atop
\overline{U}_s,\overline{V}_s,\overline{A}_s,\overline{U}^t,\overline{L}_t,\overline{\Lambda}_t,\overline{V}^t}}
2{\sum}_s\Tr(Q^TA_sP\overline{A}_s^T)+2{\sum}_t\Tr(Q^TL_t^T\overline{L}_t)\\
&\hbox{subject to}\nn
&[\overline{g};\overline{\alpha}]\in \bT,\,[\overline{h};\overline{\beta}]\in\bS,\,\overline{\mu}\geq0,\overline{\upsilon}\geq0,\,\overline{\lambda}\geq0,\,
\overline{V}^t\succeq0,\overline{S}\succeq0,\overline{T}\succeq0\nn
&\left[\begin{array}{c|c}\overline{U}_s&\overline{A}_s\cr\hline
\overline{A}_s^T&\overline{V}_s\cr\end{array}\right]\succeq0,\,\left[\begin{array}{c|c}\overline{U}^t&\overline{L}_t^T\cr\hline \overline{L}_t&\overline{\Lambda}_t\cr\end{array}\right]\succeq0\nn
&\overline{\alpha}=1,\,[\overline{g};\overline{\alpha}]\in\bS,\,\overline{\beta}=1,\,[\overline{h};\overline{\beta}]\in\bT,\,
-\overline{g}_\ell+\Tr(\overline{S}S_\ell)+\overline{\mu}_\ell=0,\,-\overline{h}_k+\Tr(\overline{T}T_k)+\overline{\upsilon}_k=0\nn
&\Tr(\overline{\Lambda}_t)-\Tr(\overline{V}_tP^TR^TR_tP)+\overline{\lambda}_t=0\nn
&\overline{U}_s=\overline{S},\,\overline{U}_t=\overline{S},\,\overline{V}_s=\overline{T},\,\overline{V}^t=\overline{T}\nonumber
\end{align}
(here and in what follows the constraints should be satisfied for all values of ``free indexes'' $s\leq S$, $t\leq T$, $\ell\leq L$, $k\leq K$). Taking into account that relation
$
\left[\begin{array}{c|c}X&Y\cr\hline Y^T&Z\cr\end{array}\right]\succeq 0$ is equivalent to $X\succeq0,Z\succeq0$, and $Y=X^{1/2}\Delta Z^{1/2}$ with $\|\Delta\|_{2,2}\leq1$, and that $[\overline{g};1]\in\bS$, $[\overline{h};1]\in\bT$ is the same as $\overline{g}\in\cS$, $\overline{h}\in\cT$,  $(D)$ boils down to
\begin{align*}
\Opt=\max\limits_{{\overline{g},\overline{h},\overline{S},\overline{T},\atop
\overline{\Delta}_s,\overline{\delta}_t,\overline{\Lambda}_t}}\bigg\{&
{\sum}_s\Tr(Q^TA_sP\overline{A}_s^T)+{\sum}_t\Tr(Q^TL_t^T\overline{L}_t):\\
&\left.\begin{array}{l}
\overline{g}\in\cT,\,\overline{h}\in\cS,\,\overline{S}\succeq0,\,\overline{T}\succeq0,\,
\Tr(\overline{S}S_\ell)\leq\overline{g}_\ell,\,\Tr(\overline{T}T_k)\leq \overline{h}_k\\
\overline{A}_s=\overline{S}^{1/2}\ov\Delta_s\overline{T}^{1/2},\,\|\ov\Delta_s\|_{2,2}\leq1,
\overline{L}_t^T=\overline{S}^{1/2}\ov\delta_t\overline{\Lambda}_t^{1/2},\,\|\ov\delta_t\|_{2,2}\leq1\\
\Tr(\overline{\Lambda}_t)\leq \Tr(\overline{S}^{1/2}P^TR^TR_tP\overline{S}^{1/2})
\end{array}\right\}
\end{align*}
or, which is the same,
\begin{align}\tag{D$'$}
\Opt=\max\limits_{{\overline{g},\overline{h},\overline{S},\overline{T}\atop
\overline{\Delta}_s,\overline{\delta}_t,\overline{\Lambda}_t,\overline{L}_t}}\bigg\{&
{\sum}_s\Tr(\overline{S}^{1/2}Q^TA_sP\overline{T}^{1/2}\ov\Delta_s^T)+2{\sum}_t\Tr(\overline{S}^{1/2}Q^TL_t^T\overline{\Lambda}_t^{1/2}\ov\delta_t^T):\\
&\left.\begin{array}{l}
\overline{g}\in\cT,\,\overline{h}\in\cS,\,\overline{S}\succeq0,\,\overline{T}\succeq0,\,\Tr(\overline{S}S_\ell)\leq\overline{g}_\ell,\,\Tr(\overline{T}T_k)\leq \overline{h}_k\\
\|\ov\Delta_s\|_{2,2}\leq1,\,\|\ov\delta_t\|_{2,2}\leq1\\
\Tr(\overline{\Lambda}_t)\leq \Tr(\overline{T}^{1/2}P^TR_t^TR_tP\overline{T}^{1/2}),\,\overline{\Lambda}_t\succeq0
\end{array}\right\}\nonumber
\end{align}
Note that for $\Delta$ and $\delta$ such that $\|\Delta\|_{2,2}\leq 1$ and $\|\delta\|_{2,2}\leq1$ one has
\[\Tr(A\Delta)\leq \|A\|_\nuc=\|\lambda(\cL[A])\|_1,\;\cL[A]=\left[\begin{array}{c|c}&\half A\cr\hline\half A^T&\cr\end{array}\right]
\]
and
\[
Tr(AB^T\delta)=\langle A,\delta^T B\rangle_\Fro\leq\|A\|_\Fro\|\delta^TB\|_\Fro\leq\|A\|_\Fro\|B\|_\Fro
\]
(here  $\|A\|_\nuc$ stands for the nuclear norm and $\lambda(A)$ for the vector of eigenvalues of a symmetric matrix $A$).
Consequently, for a feasible solution to (D$'$) it holds
$$
\Tr(\overline{S}^{1/2}Q^TA_sP\overline{T}^{1/2}\ov\Delta_s^T)\leq\|\lambda(\cL[\overline{S}^{1/2}Q^TA_sP\overline{T}^{1/2}])\|_1,\\
$$
and
$$
\Tr(\overline{S}^{1/2}Q^TL_t^T\overline{\Lambda}_t^{1/2}\ov\delta_t^T)\leq\|\overline{S}^{1/2}Q^TL_t^T\|_\Fro\|\overline{\Lambda}_t^{1/2}\|_\Fro.
$$
The latter bound combines with the last constraint in (D$'$) to imply that
$$
\Tr(\overline{S}^{1/2}Q^TL_t^T\overline{\Lambda}_t^{1/2}\ov\delta_t^T)\leq \|\overline{S}^{1/2}Q^TL_t^T\|_\Fro\|\overline{T}^{1/2}P^TR_t^T\|_\Fro,
$$
and we conclude that
\begin{align}\label{weget}
\Opt\leq \max\limits_{\overline{S},\overline{g},\overline{T},\overline{h}}\bigg\{&
{\sum}_s\left\|\lambda(\cL[\overline{S}^{1/2}Q^TA_sP\overline{T}^{1/2}])\right\|_1+
{\sum}_t\left\|\overline{S}^{1/2}Q^TL_t^T\|_\Fro\|\overline{T}^{1/2}P^TR_t^T\right\|_\Fro:\\
&\qquad\qquad\qquad\qquad\left.\begin{array}{l}
\overline{S}\succeq0,\overline{g}\in\cS,\,\Tr(\overline{S}S_\ell)\leq\overline{g}_\ell,\,\ell\leq L\\
\overline{T}\succeq0,\overline{h}\in\cT,\,\Tr(\overline{T}T_k)\leq\overline{h}_k,\,k\leq K
\end{array}\right\}\nonumber\end{align}
\paragraph{4$^o$.} We need the following result:\\
\begin{lemma}\label{lem:c1}{\cite[Lemma 2.3]{BTNMC} (cf. also \cite[Lemma 3.4.3]{BN2001})} If the ranks of all matrices $A_s$  (and thus---matrices $\overline{S}^{1/2}Q^T A_sP\overline{T}^{1/2}$) do not exceed a given $\kappa\geq1$,
then for $\omega\sim\cN(0,I_{M+N})$ one has
\[
\bE\left\{|\omega^T \cL[\overline{S}^{1/2}Q^T A_sP\overline{T}^{1/2}]\omega|\right\}\geq \|\lambda(\cL[\overline{S}^{1/2}Q^T A_sP\overline{T}^{1/2}])\|_1/\vartheta(2\kappa),
\]
with $\vartheta(\cdot)$ as described in Proposition \ref{propmain}.
\end{lemma}
Our next result is as follows (cf. \cite[Proposition B.4.12]{RO})
\begin{lemma}\label{lemfull} Let $\in\bR^{p\times q}$, $B\in\bR^{r\times q}$ and $\xi\sim\cN(Q,I_q)$. Then
$$
\bE_\xi\left\{\|A\xi\|_2\|B\xi\|_2\right\}\geq {2\over \pi}\|A\|_\Fro\|||B\|_\Fro.
$$
\end{lemma}
{\bf Proof.} Setting $A^TA=U\Diag\{\lambda\}U^T$ with orthogonal $U$ and $\zeta=U^T\xi$, we have
$$
\bE\left\{\|A\xi\|_2\|B\xi\|_2\right\}=\bE\left\{\sqrt{{{\sum}}_{i=1}^q\lambda_i[U^T\xi]_i^2}\|B\xi\|_2\right\}.
$$
The right hand side is concave in $\lambda$, so that the infimum of this function in $\lambda$ varying in the simplex ${\sum}_i\lambda_i=\Tr(A^TA)$ is attained at an extreme point. In other words, there exists vector $a\in\bR^q$ with $a^Ta=\|A\|_\Fro^2$ such that
$$
\bE\left\{\|A\xi\|_2\|B\xi\|_2\right\}\geq\bE_\xi\left\{|a^T\xi|\,\|B\xi\|_2\right\}.
$$
Applying the same argument to $\|B\xi\|_2$-factor, we can now find a vector $b\in\bR^q$, $b^Tb=\|B\|_\Fro^2$, such that
$$
\bE_\xi\left\{|a^T\xi|\,\|B\xi\|_2\right\}\geq \bE_\xi\left\{|a^T\xi|\,|b^T\xi|\right\}.
$$
It suffices to prove that the concluding quantity is $\geq 2\|a\|_2\|b\|_2/\pi$. By homogeneity, this is the same as to  prove that
if $[s;t]\sim\cN(0,I_2)$, then $\bE\{|t|\,|\cos(\phi)t+\sin(\phi)s|\}\geq {2\over\pi}$ for all $\phi\in[0,2\pi)$, which is straightforward
(for the justification, see the proof of Proposition 2.3 of \cite{ben28extended}). \qed
\par
The last building block is the following
\begin{lemma}{\cite[Lemma 6]{JuKoNe}}\label{lenlemlem}
Let
\[
\cV=\{v\in\bR^d: \exists r\in\cR: v^T R_jv\leq r_j,1\leq j\leq J\}\subset\bR^d
\]
be a basic ellitope, $W\succeq 0$ be symmetric $d\times d$ matrix such that
\[
\exists r\in\cR: \Tr(WR_j)\leq r_j,j\leq J,
\]
and $\omega\sim\cN(0,W)$. Denoting by $\rho(\cdot)$ the norm on $\bR^d$ with the unit ball $\cV$, we have
\[
\bE\{\rho(\omega)\}\leq \varkappa(J).
\]
with $\varkappa(\cdot)$ given by (\ref{varkappaJ}).
\end{lemma}

\paragraph{4$^o$} Now we can complete the proof of the second inequality in (\ref{and}). Let $\kappa\geq1$, and let $\overline{g},\overline{S},\overline{h},\overline{T}$ be feasible for the optimization problem in (\ref{weget}).  Denoting by $\|\cdot\|_\cQ$ the norm with the unit ball $\cQ$, for all $A\in\bR^{m\times n}$, $u\in\bR^m$, and $v\in\bR^n$ we have
\[
u^TAv\leq \|u\|_{\cB_*}\|Av\|_\cB\leq \|u\|_{\cB_*}\|A\|_{\cX,\cB}\|v\|_\cX,
\]
so that
for all $u\in\bR^m$ and $v\in\bR^n$
\begin{align*}
\|u\|_{\cB_*}\|v\|_\cX\|\cA\|_{\cX,\cB}&\geq
\max\limits_{{\epsilon_s,|\epsilon_s|\leq1,
\atop\delta_t,\|\delta_t\|_{2,2}\leq1}}\left[{\sum}_s\epsilon_su^TA_sv+{\sum}_tu^TL_t^T\delta_tR_tv\right]
\\&={\sum}_s|u^TA_sv|+{\sum}_t\|L_tu\|_2\|R_tv\|_2.
\end{align*}
Thus, for all $\overline{g},\overline{S},\overline{h},\overline{T}$ which are feasible for \rf{weget} and $\xi\in\bR^M,\,\eta\in\bR^N$,
\begin{align}
\|\overline{S}^{1/2}\xi\|_\cZ\|\overline{T}^{1/2}\eta\|_\cY\|\cA\|_{\cX,\cB}&\geq
\|Q\overline{S}^{1/2}\xi\|_{\cB_*}
\|P\overline{T}^{1/2}\eta\|_\cX \|\cA\|_{\cX,\cB}\hbox{\ [due to $\cB_*=Q\cZ,\cX=P\cY$]}\nn
&\geq {\sum}_s|\xi^T\overline{S}^{1/2}Q^TA_sP\overline{T}^{1/2}\eta|+{\sum}_t\|L_tQ\overline{S}^{1/2}\xi\|_2\|R_tP\overline{T}^{1/2}\eta\|_2\nn
&={\sum}_s|[\xi;\eta]^T\cL[\overline{S}^{1/2}Q^TA_sP\overline{T}^{1/2}][\xi;\eta]|\nn
&+{\sum}_t\|[L_tQ\overline{S}^{1/2},0_{p_t\times N}][\xi;\eta]\|_2
\|[0_{q_t\times M},R_tP\overline{T}^{1/2}][\xi;\eta]\|_2.
\label{eq*}
\end{align}
As a result, for $[\xi;\eta]\sim\cN(0,I_{M+N})$, applying the bounds of Lemmas \ref{lem:c1} and \ref{lemfull},
\begin{align*}
&\bE\left\{\left\|\overline{S}^{1/2}\xi\right\|_\cZ\right\}\bE\left\{\left\|\overline{T}^{1/2}\eta\right\|_\cY\right\}\|\cA\|_{\cX,\cB}
=\bE\left\{\left\|\overline{S}^{1/2}\xi\right\|_\cZ\|\overline{T}^{1/2}\eta\|_\cY \|\cA\|_{\cX,\cB}\right\}\\
&\geq {\sum}_s\bE\left\{\left|[\xi;\eta]^T\cL[\overline{S}^{1/2}Q^TA_sP\overline{T}^{1/2}][\xi;\eta]\right|\right\}\\
&\quad+{\sum}_t\bE\left\{\left\|[L_tQ\overline{S}^{1/2},0_{p_t\times N}][\xi;\eta]\right\|_2
\left\|[0_{q_t\times M},R_tP\overline{T}^{1/2}][\xi;\eta]\right\|_2\right\}
\\&\geq \vartheta(2\kappa)^{-1}\sum_s\left\|\lambda\left(\cL[\overline{S}^{1/2}Q^TA_sP\overline{T}^{1/2}]\right)\right\|_1
+\tfrac{2}{\pi}\sum_t\|L_tQ\overline{S}^{1/2}\|_\Fro\|R_tP\overline{T}^{1/2}\|_\Fro.
\end{align*}
Besides this, by Lemma \ref{lenlemlem} we have
\[
\bE\left\{\|\overline{S}^{1/2}\xi\|_\cZ\right\}\leq\varkappa(L),\;\;
\bE
\left\{\|\overline{T}^{1/2}\eta\|_\cY\right\}\leq\varkappa(K)
\]
due to the fact that $\overline{g},\overline{S},\overline{h}$ and $\overline{T}$ are feasible for (\ref{weget}). This combines with \rf{eq*} to imply that
the value $\varkappa(L)\varkappa(K)\|\cA\|_{\cX,\cB}$ is lower bounded with the quantity
\[
\max\left[\vartheta(2\kappa),\pi/2\right]^{-1}
\left[{\sum}_s\left\|\lambda\left(\cL[\overline{S}^{1/2}Q^TS_sP\overline{T}^{1/2}]\right)\right\|_1
+{\sum}_t\|\overline{S}^{1/2}Q^TL_t^T\|_\Fro\|\overline{T}^{1/2}P^TR_t^T\|_\Fro\right].
\]
Invoking the inequality in (\ref{weget}), we arrive at the second inequality in (\ref{and}). The above reasoning assumed that $\kappa\geq1$, with evident simplifications,
it is applicable to the case of $\kappa=0$ as well. \qed
\subsection{Proof of Proposition \ref{prop912}}
We put $S=q_\s$ and $T=q-q_\s$. In the situation of Proposition \ref{prop912} we want to tightly upper-bound quantity
$$
\begin{array}{rcl}
\ms(H)&=&\max\limits_{x\in \cX,\eta\in\cU}\left\|H^TD[\eta]x\right\|\\
&=&\max\limits_{\ell\leq L}\max\limits_{x\in\cX,\eta\in\cU}\left\{\sqrt{[H^TD[\eta]x]^TR_\ell [H^TD[\eta]x]}\right\}\\
&=&\max\limits_{\ell\leq L}\|\cA_\ell[H]\|_{\cX,2},\\
\end{array}
$$
where $\|\cdot\|_{\cX,2}$ is the operator norm induced by $\|\cdot\|_\cX$ on the argument and $\|\cdot\|_2$ on the image space and the uncertain matrix $\cA_\ell[H]$ is given by
$$
\begin{array}{c}
\cA_\ell=\bigg\{{\sum}_{s=1}^{S}\delta_s
\underbrace{R_\ell^{1/2}H^TP_s^TQ_s}_{=:A_{s\ell}[H]}
+{\sum}_{t=1}^{T}\underbrace{R_\ell^{1/2}H^TP_{S+t}^T}_{L_{t\ell}^T[H]}\Delta_s \underbrace{Q_{S+t}}_{=:R_t}:
\\
\multicolumn{1}{r}{\begin{array}{ll}
|\delta_s|\leq1&,1\leq s\leq S\\
\|\Delta_s\|_{2,2}\leq 1&,1\leq t\leq T\\
\end{array}\bigg\}}
\end{array}
$$
It follows that
$$
\ms(H)=\max\limits_{\ell\leq L} \|\cA_\ell[H]\|_{\cX,2},
$$
and Proposition \ref{propmain} provides us with the efficiently computable convex in $H$ upper bound $\ov\ms(H)$ on $\ms(H)$:
\begin{align*}
\ov\ms(H)&=\max\limits_{\ell\leq L}\Opt_\ell(H),\\
\Opt_\ell(H)&=\min\limits_{\mu,\upsilon,\lambda,U_s,V_s,U^t,V^t}
\Big\{\half[\mu+\phi_\cT(\upsilon)]:\,\mu\geq0,\,\upsilon\geq0,\,\lambda\geq0\\
&\qquad\qquad\qquad\left.\begin{array}{l}
\left[\begin{array}{c|c}U_s&-A_{s\ell}[H]P\cr\hline -P^TA_{s\ell}^T[H]&V_s\cr\end{array}\right]\succeq0\\
\left[\begin{array}{c|c}U^t&-L_{t\ell}^T[H]\cr\hline -L_{t\ell}[H]&\lambda_tI_{p_{q_\s+t}}\cr\end{array}\right]\succeq0\\
V^t-\lambda_tP^TR_t^TR_tP\succeq0\\
\mu I_\nu-{\sum}_sU_s-{\sum}_tU^t\succeq0\\
{\sum}_k\upsilon_k T_k-{\sum}_sV_s-{\sum}_tV^t\succeq0
\end{array}\right\}
\end{align*}
and tightness factor of this bound does not exceed $\max[\vartheta(2\kappa),\pi/2]$ where $\kappa=\max\limits_{\alpha\leq q_\s}\min[p_\alpha,q_\alpha]$.\qed
\subsection{Spectratopic version of  Proposition \ref{propmain}}
Proposition \ref{propmain} admits a ``spectratopic version,'' in which ellitopes $\cX$ and $\cB_*$ given by (\ref{ellitopes}) are replaced by the pair of {\em spectratopes}
\begin{subequations} \label{spectratopes}
\begin{align}
&\begin{array}{rcl}\cX&=&\{Py: y\in\cY\}\subset\bR^n, \cY=\{y\in\bR^N\ \&\ \exists t\in\cT: \,T_k[y]^2\preceq t_kI_{f_k},k\leq K\},\\
&&T_k[y]=\sum_{j=1}^N y_jT^{jk},\,T^{jk}\in\bS^{f_k},\sum_kT_k^2[y]\succ0\,\,\forall y\neq0\end{array}
\label{spectratopes.a}\\
&\begin{array}{rcl}
\cB_*&=&\{Qz:z\in\cZ\}\subset\bR^m,\,\cZ=\{z\in\bR^M:\exists s\in\cS:\, S_\ell^2[z]\preceq s_\ell I_{d_\ell},\,\ell\leq L\},\\
&&S_\ell[z]=\sum_{j=1}^M z_jS^{jk\ell},\,S^{jk\ell}\in\bS^{d_\ell},\sum_\ell S_\ell^2[z]\succ0\,\,\forall z\neq0 \end{array}\label{spectratopes.b}
\end{align}
\end{subequations}
The spectratopic version of the statement reads as follows:
\begin{proposition}\label{propmainspectr} Given uncertain matrix (\ref{cA}) and spectratopes \rf{spectratopes.a} and \rf{spectratopes.b}, consider convex optimization problem
\begin{subequations}\label{cprbspectr}
\begin{align}
\Opt&=\min\limits_{\mu,\upsilon,\lambda,U_s,V_s,U^t,V^t}\bigg\{\half[\phi_\cS(\lambda[\mu])+\phi_\cT(\lambda[\upsilon])]:
\nn
&\hbox{subject to}\nn
&\mu=\{M_\ell\in\bS^{d_\ell}_+,\ell\leq L\},\,\upsilon=\{\Upsilon_k\in\bS^{f_k}_+,k\leq K\},\,\lambda\geq0\nn
&\left[\begin{array}{c|c}U_s&-Q^TA_sP\cr\hline -P^TA_s^TQ&V_s\cr\end{array}\right]\succeq0\label{cprbspectr.a}\\
&\left[\begin{array}{c|c}U^t&-Q^TL_t^T\cr\hline -L_tQ&\lambda_tI_{p_t}\cr\end{array}\right]\succeq0,\,V^t-\lambda_tP^TR_t^TR_tP\succeq0\label{cprbspectr.b}\\
&\sum_\ell S^{+,*}_\ell[M_\ell]-\sum_sU_s-\sum_tU^t\succeq0\label{cprbspectr.c}\\
&\sum_k T_k^{+,*}[\Upsilon_k]-\sum_sV_s-\sum_tV^t\succeq0\label{cprbspectr.d}
\end{align}
\end{subequations}
where
$$
\lambda[\zeta]=[\Tr(Z_1);...;\Tr(Z_I)]\hbox{\ for\ }\zeta=\{Z_i\in\bS^{k_i},i\leq I\}
$$
and
$$ S^{+,*}_\ell[V]=\left[\Tr(V S^{i\ell}S^{j\ell})\right]_{i,j \leq M}\hbox{\ for\ }V\in\bS^{d_\ell},\,
T^{+,*}_k[U]=\left[\Tr(U T^{ik}T^{jk})\right]_{i,j \leq N}\hbox{\ for\ }U\in\bS^{f_k}.$$
\par\noindent
Problem \rf{cprbspectr} is strictly feasible and solvable, and
\[
\|\cA\|_{\cX,\cB}\leq\Opt\leq \varsigma\left({\sum}_kf_k\right)\varsigma\left({\sum}_\ell d_\ell\right)
\max\left[\vartheta(2\kappa),\pi/2\right]\|\cA\|_{\cX,\cB}
\]
where
$\vartheta$ and $\kappa$ are the same as in Proposition \ref{propmain}
and
\[
\varsigma(J)=2\sqrt{2\ln(2J)}.
\]
\end{proposition}
{\bf Proof.} For $Y\in\bS^M$ and $X\in\bS^N$ let us set
$$
S_\ell^+[Y]=\sum_{i,j=1}^MY_{ij}S^{i\ell}S^{j\ell},\,\,T_k^+[X]=\sum_{i,j=1}^NX_{ij}T^{ik}T^{jk},
$$
so that
\begin{equation}\label{ident1}
S_\ell^+[zz^T]=S_\ell^2[z],\,T_k^+[yy^T]=T_k^2[y]
\end{equation}
and
\begin{equation}\label{ident2}
\begin{array}{l}
\Tr(VS_\ell^+[Y])=\Tr(S_\ell^{+,*}[V]Y)\hbox{\ for\ }V\in\bS^{d_\ell},Y\in\bR^M,\\
\Tr(UT_k^+[X])=\Tr(T_k^{+,*}[U]X)\hbox{\ for\ }U\in\bS^{f_k},X\in\bR^N.
\end{array}\end{equation}
The proof of Proposition \ref{propmainspectr} is obtained from that (below referred to as ``the proof'') of Proposition \ref{propmain} by the following modifications:
\begin{enumerate}
\item All references to (\ref{cprb}) should be replaced with references to (\ref{cprbspectr}). Item [b] in 1$^o$ of the proof now reads
\begin{quotation} [b$'$] $x=Py$ with $T_k^2[y]\preceq \tau_kI_{f_k}$, $k\leq K$, for some $\tau\in\cT$ and $u=Qz$ for some $z$ such that $S_\ell^2[z]\preceq \varsigma_\ell I_{d_\ell}$, $\ell\leq L$, for some $\varsigma\in\cS$.
\end{quotation}
The last three lines in the chain \rf{chain74} are replaced with
\begin{align*}
\gamma&\leq {1\over 2}\left[\sum_\ell\Tr([zz^T] S_\ell^{+,*}[M_\ell])+\sum_k\Tr([yy^T]T_k^{+,*}[\Upsilon_k])\right]\hbox{\ [by \rf{cprbspectr.c} and\rf{cprbspectr.d}]}\\
&={1\over 2}\left[\sum_\ell\Tr(S_\ell^2[z]M_\ell)+\sum_k\Tr(T_k^2[y]\Upsilon_k)\right]\hbox{\ [by (\ref{ident1}) and (\ref{ident2})]}\\
&\leq  {1\over 2}\left[\sum_\ell\varsigma_\ell\Tr(M_\ell)+\sum_k\tau_k\Tr(\Upsilon_k)\right]\hbox{\ [due to (b$'$) and $M_\ell\succeq0,\Upsilon_k\succeq0$]}\\
&\leq \half[\phi_\cS(\lambda[\mu])+\phi_\cT(\lambda[\upsilon])]\hbox{\ [since $\varsigma\in\cS,\tau\in\cT$]}.
\end{align*}
\item Constraints \rf{tagP.ef} in (P) now read
$$
\left[\sum_\ell S_\ell^{+,*}[M_\ell]-\sum_sU_s-\sum_tU^t\right]^{\hbox{\scriptsize$\overline{S}$}}\succeq0,\;\;
\left[\sum_kT_k^{+,*}[\Upsilon_k]-\sum_sV_s-\sum_tV^t\right]^{\hbox{\scriptsize$\overline{T}$}}\succeq0.
$$
As a result, (\ref{weget}) becomes
\begin{align}
\label{wegetspectr}
\Opt\leq \max\limits_{\overline{S},\overline{g},\overline{T},\overline{h}}\bigg\{&
\sum_s\left\|\lambda(\cL[\overline{S}^{1/2}Q^TA_sP\overline{T}^{1/2}])\right\|_1
+
\sum_t\|\overline{S}^{1/2}Q^TL_t^T\|_\Fro\|\overline{T}^{1/2}P^TR_t^T\|_\Fro:\nn
&\left.\qquad\qquad\qquad\qquad\qquad\begin{array}{l}
\overline{S}\succeq0,\overline{g}\in\cS,S_\ell^+[\overline{S}]\preceq\overline{g}_\ell I_{d_\ell},\,\ell\leq L\\
\overline{T}\succeq0,\overline{h}\in\cT,T_k^+[\overline{T}]\preceq\overline{h}_k I_{f_k},\,k\leq K
\end{array}\right\}
\end{align}
\item The role of Lemma \ref{lenlemlem} in the proof is now played by the following fact.
\begin{lemma}\label{sslenlemlem}{\cite[Lemma 8]{JuKoNe}}
Let
\[
\cV=\{v\in\bR^d: \exists r\in\cR: R_j^2[v]\preceq r_jI_{\nu_j},1\leq j\leq J\}\subset\bR^d
\]
be a basic spectratope, $W\succeq 0$ be symmetric $d\times d$ matrix such that
\[
\exists r\in\cR: R^+_j[W]\preceq  r_j I_{\nu_j},j\leq J,
\]
and $\omega\sim\cN(0,W)$. Denoting by $\rho(\cdot)$ the norm on $\bR^d$ with the unit ball $\cV$, we have
\[
\bE\{\rho(\omega)\}\leq \varsigma\left({\sum}_j\nu_j\right),\;
\varsigma(F)=2\sqrt{2\ln(2F)}.\\
\]
\end{lemma}
\end{enumerate}



\end{document}